\documentclass[10pt,a4paper,twoside]{article}

\usepackage{amsmath,amsfonts,amssymb,amscd,graphics}
\usepackage{theorem}
\usepackage[english]{babel}

\theoremstyle{plain}
\newtheorem{theo}{Theorem}
\newtheorem{definition}[theo]{Definition}

\newtheorem{prop}[theo]{Proposition}
\newtheorem{lemme}[theo]{Lemma}
\newtheorem{cor}[theo]{Corollary}

{\theorembodyfont{\upshape}
}

{\theorembodyfont{\upshape}
}

{\theorembodyfont{\upshape}
\newtheorem{rem}[theo]{Remark}}

\def\demo{\noindent \textsc{Proof} \\}
\def\findem{\hfill$\square$\vskip 13pt}

\title{Stability of solutions of BSDEs with random terminal time}
\author{Sandrine TOLDO\footnote{{\it Email address :} {\rm sandrine.toldo@math.univ-rennes1.fr}} \\
{\it IRMAR, Universit\'e Rennes 1, Campus de Beaulieu, 35042 Rennes Cedex, France}}
\date{}
\begin{document}

\maketitle

\begin{abstract}  
In this paper, we study the stability of the solutions of Backward Stochastic Differential Equations (BSDE for short) with an almost surely
finite random terminal time. More precisely, we are going to show that if $(W^n)$ is a sequence of scaled random walks or a sequence of 
martingales that converges to a Brownian motion $W$ and if $(\tau^n)$ is a sequence of stopping times that converges to a stopping time 
$\tau$, then the solution of the BSDE driven by $W^n$ with random terminal time $\tau^n$ converges to the solution of the BSDE driven by 
$W$ with random terminal time $\tau$.
\end{abstract}
\noindent
\textbf{Keywords :} Backward Stochastic Differential Equations (BSDE), Stability of BSDEs, Weak convergence of filtrations, Stopping times 
\maketitle
\section*{Introduction}
We want to make an approximation of the solutions of a backward stochastic differential equation (BSDE for short) with an almost surely
finite random terminal time $\tau$ like 
$$Y_{t \wedge \tau} = \xi + \int_{t \wedge \tau}^{\tau} f(s,Y_s,Z_s)ds - \int_{t \wedge \tau}^{\tau} Z_s dW_s, \ t \geqslant 0. $$
\indent The robustness of numerical methods to approximate solutions of BSDEs had already been studied in the case of a deterministic
terminal time. 
Antonelli and Kohatsu-Higa in \cite{AK2000} and also Coquet, Mackevi{\v{c}}ius and M\'emin in \cite{CMM98} and \cite{CMM98cor} proposed 
approximation schemes that use discretization of filtrations and convergence of filtrations. 
Briand, Delyon and M\'emin in \cite{BDM_Donsker} and Ma, Protter, San Mart{\'{\i}}n and Torres in \cite{MPSMT} have approximated the 
Brownian motion by a scaled random walk. Then, in \cite{BDM_mart}, Briand, Delyon and M\'emin have studied another case~: they approach 
the Brownian motion by a sequence of martingales. \\
\indent We are interested in BSDEs with random terminal time because there have strong links with Partial Differential Equations as it is
explained by Peng in \cite{Peng91}. As for BSDEs with deterministic terminal time, we study the robustness of numerical methods to
approximate the solutions of those BSDEs. In this paper, we shall approximate the Brownian motion either by a scaled random walk, either by
a sequence of martingales. In this study, we need moreover to approximate the random almost surely finite terminal time $\tau$ by a 
sequence of stopping times $(\tau^n)_n$. \\
\indent In Section \ref{stab_donsker}, we approximate the Brownian motion by a scaled random walk. First, we shall state the
problem and study the properties of existence and uniqueness of the solutions of the BSDEs. Then, we will deal with the convergence of the
solutions. To end this part, we give an example of the convergence result for hitting times. \\
\indent In Section \ref{stab_mart}, we approximate the Brownian motion by a sequence of martingales. We shall see some
generalizations of the results of Section \ref{stab_donsker} : existence and uniqueness of the solutions of the BSDEs under study, 
convergence of the solutions. Moreover, we will illustrate these results by the case of discretizations of a Brownian motion and hitting 
times. \\
\indent For technical reasons, we need some results about convergence of stopped filtrations. So, in Appendix \ref{cv_stop_filtr}, we will 
deal with stopped filtrations and stopped processes. We are going to establish a link between the convergence of a sequence of stopped 
processes and the convergence of the associated stopped filtrations. \\

In what follows, we are given a probability space $(\Omega, \mathcal{A}, \mathbb{P})$. 
Unless otherwise specified, every $\sigma$-field will be supposed to be included in $\mathcal{A}$, every process will be indexed by $\mathbb{R}^+$ 
and taking values in $\mathbb{R}$, every filtration will be indexed by $\mathbb{R}^+$. $\mathbb{D}=\mathbb{D}(\mathbb{R}^+)$ denotes the space of c\`adl\`ag 
functions from $\mathbb{R}^+$ to $\mathbb{R}$. We endow $\mathbb{D}$ with the Skorokhod topology. If $X$ is a process and $\tau$ a stopping time, we denote 
by $X^\tau$ the corresponding stopped process, $i.e.$~for every $t$, $X^\tau_t = X_{t \wedge \tau}$. 
\begin{definition}
Let $(\mathcal{F}_t)_{t \geqslant 0}$ be a filtration and $\tau$ a $\mathcal{F}$-stopping time. 
We define the $\sigma$-field $\mathcal{F}_\tau$ by 
$$\mathcal{F}_\tau = \{A \in \mathcal{F}_\infty : A \cap \{\tau \leqslant s \} \in \mathcal{F}_s, \forall s \} \text{~~where~~} \mathcal{F}_\infty = \bigvee_t \mathcal{F}_t$$ 
and the stopped filtration $\mathcal{F}^\tau$ by 
$$\mathcal{F}^\tau_t = \mathcal{F}_{\tau \wedge t} = \{A \in \mathcal{F}_t : A \cap \{\tau \leqslant s \} \in \mathcal{F}_s, \forall s \leqslant t\},$$
for every $t$.
\end{definition}
\indent For technical background about Skorokhod topology, the reader may refer to Billingsley \cite{Bill} or Jacod and Shiryaev 
\cite{JS}.\\


\section{Stability of BSDEs when the Brownian motion is approximated by a scaled random walk}
\label{stab_donsker}
\label{hb2.1}


\subsection{Statement of the problem}

Let $f:\mathbb{R}^2 \to \mathbb{R}$ be a Lipschitz function : there exists $K \in \mathbb{R}^+$ such that, for every $(y,z), (y',z') \in \mathbb{R}^2$, we have~:
$$|f(y,z) - f(y',z')| \leqslant K[|y-y'|+|z-z'|].$$
For clarity's sake, we consider a time-independent generator $f$ but the results of Section \ref{stab_donsker} remain true if $f$ is
time-dependent as it is explained in Remark \ref{rem_f}. \\

We also suppose that $f$ is bounded and that $f$ fills the following property of monotonicity w.r.t. $y$ : there exists $\mu > 0$ such that 
$$\forall (y,z),(y',z) \in \mathbb{R}^+ \times \mathbb{R}^2, (y-y')(f(y,z)-f(y',z)) \leqslant - \mu (y-y')^2.$$

Let $W$ be a Brownian motion and $\mathcal{F}$ its natural filtration. Let $\tau$ be a $\mathcal{F}$-stopping time almost surely finite. \\
We consider the following stochastic differential equation :
\begin{equation}
\label{Y}
Y_{t \wedge \tau} = \xi + \int_{t \wedge \tau}^{\tau} f(Y_s,Z_s)ds - \int_{t \wedge \tau}^{\tau} Z_s dW_s, \ t \geqslant 0, 
\end{equation}
where $\xi$ is a bounded $\mathcal{F}_\tau$-mesurable random variable. \\

\begin{definition}
\label{hb_def_sol}
We call a solution of the BSDE (\ref{Y}) a pair $(Y,Z)$ of progressively measurable processes verifying the equation (\ref{Y}) such that
$Y_t=\xi$ and $Z_t=0$ on the set $\{t>\tau\}$ and :
$$\mathbb{E} \left[ \sup_{t \in \mathbb{R}^+} e^{-2\mu t} |Y_t|^2 \right] < + \infty \text{~and~} 
	\forall t, \  \mathbb{E} \left[ \int_0^{t \wedge \tau} |Z_t|^2 dt \right] < + \infty.$$
\end{definition}

According to Theorem 2.1 of Royer in \cite{article_Manuela}, the BSDE (\ref{Y}) has a unique pair solution $(Y,Z)$ in the set of 
processes such that $Y$ is continuous and uniformly bounded. \\

In this section, we approximate equation (\ref{Y}) on the following way. 
We consider the sequence of scaled random walks $(W^n)_{n \geqslant 1}$ defined by :
$$W^n_t=\frac{1}{\sqrt{n}}\sum_{k=1}^{[nt]} \varepsilon_k^n, \ t \geqslant 0$$
where $(\varepsilon_k^n)_{k \in \mathbb{N}^*}$ is a sequence of i.i.d. symmetric Bernoulli variables. Let $\mathcal{F}^n$ be the natural
filtrations of $W^n$, $n \geqslant 1$. We have $\mathcal{F}^n_t = \sigma(\varepsilon^n_k, k \leqslant [nt])$. 
Let $(\tau^n)_n$ be a sequence of bounded $(\mathcal{F}^n)$-stopping times. 
For each $n$, we can find $T_n \in \mathbb{N}$ such that $\tau^n \leqslant T_n$. \\

Then, for each $n$, we consider the following equation :
\begin{equation}
\label{Yn*}
\begin{array}{l}
y^n_{\frac{k}{n} \wedge \tau^n}=y^n_{\frac{k+1}{n} \wedge \tau^n} + \frac{1}{n} \mathbf{1}_{\{\tau^n \geqslant \frac{k}{n}\}}
f \left(y^n_{\frac{k}{n} \wedge \tau^n},z^n_{\frac{k+1}{n} \wedge \tau^n} \right) - z^n_{\frac{k+1}{n} \wedge \tau^n} \frac{1}{\sqrt{n}}
\varepsilon^n_{k+1},\ k=0, \ldots, (n-1)T_n \\
y^n_{\tau^n}=\xi^n
\end{array}
\end{equation}
where $(\xi^n)$ is a sequence of $(\mathcal{F}^n_{\tau^n})$-measurable integrable random variables.\\

By a solution of equation (\ref{Yn*}), we mean a discrete process $\{y^n_{\frac{k}{n}},z^n_{\frac{k+1}{n}}\}_{k \geqslant 0}$ 
that satisfies (\ref{Yn*}), such that $y^n_{\frac{k}{n}}=\xi^n$ and $z^n_{\frac{k}{n}}=0$ on the set $\{\tau^n < \frac{k}{n}\}$ and such
that $\{y^n_{\frac{k}{n} \wedge \tau^n},z^n_{\frac{k+1}{n} \wedge \tau^n}\}_{k \geqslant 0}$ is $\mathcal{F}^{n,\tau^n}$-adapted. \\

\begin{prop}
\label{exist_unic_Yn*}
Equation (\ref{Yn*}) has a unique solution $(y^n,z^n)$. 
\end{prop}

\demo
We are going to build this solution on a converse iterative way. \\
For $k=nT_n$, let us put $y^n_{\frac{k}{n} \wedge \tau^n} = \xi^n$ and $z^n_{\frac{k+1}{n} \wedge \tau^n} = 0$. \\
Let us suppose that, for a given $k$, we have built 
$\left(y^n_{\frac{k+1}{n} \wedge \tau^n},z^n_{\frac{k+2}{n} \wedge \tau^n}\right)$. 
Using equation (\ref{Yn*}), we are going to determinate 
$\left(y^n_{\frac{k}{n} \wedge \tau^n},z^n_{\frac{k+1}{n} \wedge \tau^n}\right)$. \\

Let us begin by giving an expression of $z^n_{\frac{k+1}{n} \wedge \tau^n}$.\\
Multiplying equation (\ref{Yn*}) by $\sqrt{n} \varepsilon^n_{k+1} \mathbf{1}_{\{\tau^n \geqslant \frac{k}{n}\}}$ and taking the conditional 
expectation with respect to $\mathcal{F}^{n,\tau^n}_{k/n}$, we have :
\begin{eqnarray*}
\lefteqn{z^n_{\frac{k+1}{n} \wedge \tau^n}\mathbf{1}_{\{\tau^n \geqslant \frac{k}{n}\}}} \\
& = & \mathbb{E}\left[z^n_{\frac{k+1}{n} \wedge \tau^n}\mathbf{1}_{\{\tau^n \geqslant \frac{k}{n}\}} \big| \mathcal{F}^{n,\tau^n}_{\frac{k}{n}}\right] \\
& = & 
\sqrt{n} \mathbb{E}\left[y^n_{\frac{k+1}{n} \wedge \tau^n}\varepsilon^n_{k+1}\mathbf{1}_{\{\tau^n \geqslant \frac{k}{n}\}} 
			\big| \mathcal{F}^{n,\tau^n}_{\frac{k}{n}}\right] 
+ \frac{1}{\sqrt{n}} \mathbb{E}\left[ {\bf 1}_{\{\tau^n \geqslant \frac{k}{n}\}}
	f \left(y^n_{\frac{k}{n} \wedge \tau^n},z^n_{\frac{k+1}{n} \wedge \tau^n} \right)\varepsilon^n_{k+1} \big|
	\mathcal{F}^{n,\tau^n}_{\frac{k}{n}}\right] \\
&& \quad \quad - 
\sqrt{n} \mathbb{E}\left[y^n_{\frac{k}{n} \wedge \tau^n}\varepsilon^n_{k+1}\mathbf{1}_{\{\tau^n \geqslant \frac{k}{n}\}} 
			\big| \mathcal{F}^{n,\tau^n}_{\frac{k}{n}}\right] \\
& = & \sqrt{n} \mathbb{E}\left[y^n_{\frac{k+1}{n} \wedge \tau^n}\varepsilon^n_{k+1} \big| \mathcal{F}^{n,\tau^n}_{\frac{k}{n}}\right]
	\mathbf{1}_{\{\tau^n \geqslant \frac{k}{n}\}}
\end{eqnarray*}
because $y^n_{\frac{k}{n} \wedge \tau^n}$ and $z^n_{\frac{k+1}{n} \wedge \tau^n}$ must be $\mathcal{F}^{n,\tau^n}_{\frac{k}{n}}$-measurable, 
$\tau^n$ is a $\mathcal{F}^n$-stopping time and $\varepsilon^n_{k+1}$ is independent from $\mathcal{F}^{n,\tau^n}_{\frac{k}{n}}$ and is centered. \\
Moreover, by definition of a solution, for every $k$, 
$z^n_{\frac{k+1}{n} \wedge \tau^n} = z^n_{\frac{k+1}{n} \wedge \tau^n} \mathbf{1}_{\{\tau^n \geqslant \frac{k}{n}\}}$. 
So, we put :
$$z^n_{\frac{k+1}{n} \wedge \tau^n}
	=\sqrt{n}\mathbb{E}\left[y^n_{\frac{k+1}{n} \wedge \tau^n}\varepsilon^n_{k+1} \Big| \mathcal{F}^n_{\frac{k}{n}} \right]
	\mathbf{1}_{\{\tau^n \geqslant \frac{k}{n}\}}.$$
We point out that $z^n_{\frac{k+1}{n} \wedge \tau^n}$ is $\mathcal{F}^{n,\tau^n}_{\frac{k}{n}}$-measurable. \\
	
Now, we have to determinate $y^n_{\frac{k}{n} \wedge \tau^n}$.\\
On the set $\{ \tau^n < \frac{k}{n}\}$, 
	$y^n_{\frac{k}{n} \wedge \tau^n} = y^n_{\frac{k+1}{n} \wedge \tau^n} = y^n_{\tau^n} = \xi^n$. \\
On $\{ \tau^n \geqslant \frac{k}{n}\}$, 
$y^n_{\frac{k}{n} \wedge \tau^n}=y^n_{\frac{k+1}{n} \wedge \tau^n} + \frac{1}{n} 
f \left(y^n_{\frac{k}{n} \wedge \tau^n},z^n_{\frac{k+1}{n} \wedge \tau^n} \right) - z^n_{\frac{k+1}{n} \wedge \tau^n} \frac{1}{\sqrt{n}}
\varepsilon^n_{k+1} $ and we can write it  
$$y^n_{\frac{k}{n} \wedge \tau^n} = \varphi \left(y^n_{\frac{k}{n} \wedge \tau^n}\right)$$ 
with $\varphi(y)=y^n_{\frac{k+1}{n} \wedge \tau^n} + \frac{1}{n} 
f \left(y,z^n_{\frac{k+1}{n} \wedge \tau^n} \right) - z^n_{\frac{k+1}{n} \wedge \tau^n} \frac{1}{\sqrt{n}} \varepsilon^n_{k+1} $. 
As $f$ is $K$-Lipschitz in $y$, we have, for every $y,y'$ :
$$|\varphi(y)-\varphi(y')| \leqslant \frac{K}{n}|y-y'|.$$
So, for $n$ large enough (notice that this range does not depend on $k$), $\frac{K}{n}<1$ and $\varphi$ is a contraction.
Then, the equation $y^n_{\frac{k}{n} \wedge \tau^n} = \varphi \left(y^n_{\frac{k}{n} \wedge \tau^n}\right)$ has a 
unique solution for $n$ large enough, according to a fix point Theorem. By construction, $y^n_{\frac{k}{n} \wedge \tau^n}$ is
$\mathcal{F}^{n,\tau^n}_{\frac{k+1}{n}}$-measurable. But, using the predictable representation property, we have 
$y^n_{\frac{k+1}{n} \wedge \tau^n}- z^n_{\frac{k+1}{n} \wedge \tau^n} \frac{1}{\sqrt{n}} \varepsilon^n_{k+1} 
= \mathbb{E}[y^n_{\frac{k+1}{n} \wedge \tau^n}|\mathcal{F}^{n,\tau^n}_{\frac{k}{n}}]$. So $y^n_{\frac{k}{n} \wedge \tau^n}$ is independant from
$\varepsilon^n_{k+1}$ and $y^n_{\frac{k}{n} \wedge \tau^n}$ is $\mathcal{F}^{n,\tau^n}_{\frac{k}{n}}$-measurable. \\
Hence, the equation (\ref{Yn*}) has a unique solution. 
\findem
 
Now, we define continuous time processes $Y^n$ and $Z^n$ by $Y^n_t=y^n_{\frac{[nt]}{n} \wedge \tau^n}$, 
$Z^n_t=z^n_{\frac{\lfloor nt \rfloor}{n} \wedge \tau^n}$, for every $t \in \mathbb{R}^+$ where $\lfloor x \rfloor=(x-1)^+$ if $x$ is an 
integer, $[x]$ otherwise. The processes $Y^n$ and $Z^n$ are constant on the 
intervals $[k/n,(k+1)/n[$ and $]k/n,(k+1)/n]$ respectively and satisfy the following equation :
\begin{equation}
\label{Yn}
Y^n_{t}= \xi^n + \int_{t \wedge \tau^n}^{\tau^n} f(Y^n_{(s \wedge \tau^n)-},Z^n_{s \wedge \tau^n})dA^n_s 
- \int_{t \wedge \tau^n}^{\tau^n} Z^n_{s \wedge \tau^n} dW^n_s, 
\end{equation}
where $A^n_s=\frac{[ns]}{n}$.\\

\begin{rem}
\label{rem_f}
If $f$ is time-dependent such that for every $(y,z)$, we suppose that $\{f(t,y,z)\}_{t \geqslant 0}$ is progressively measurable, bounded,
$K$-Lipschitz in $y$ and $z$ and verify the condition of monotonicity w.r.t. $y$ given before. Under these assumptions, all the results of 
Section \ref{stab_donsker} remain true and the proofs are the same. In that case, Equation (\ref{Yn*}) becomes :
\begin{eqnarray*}
y^n_{\frac{k}{n} \wedge \tau^n} & = &y^n_{\frac{k+1}{n} \wedge \tau^n} + \frac{1}{n} \mathbf{1}_{\{\tau^n \geqslant \frac{k}{n}\}}
f \left(\frac{k}{n} \wedge \tau^n,y^n_{\frac{k}{n} \wedge \tau^n},z^n_{\frac{k+1}{n} \wedge \tau^n} \right) 
- z^n_{\frac{k+1}{n} \wedge \tau^n} \frac{1}{\sqrt{n}}
\varepsilon^n_{k+1},\ k=0, \ldots, (n-1)T_n \\
y^n_{\tau^n} & = & \xi^n.
\end{eqnarray*}
\end{rem}


\subsection{Convergence of the solutions}

The aim of this section is to prove the following result of convergence of the solutions :
\begin{theo}
\label{theo_Donsker}
Let $(Y,Z)$ be the solution of the BSDE (\ref{Y}) and $(Y^n,Z^n)$ be the processes constant on the 
intervals $[k/n,(k+1)/n[$ and $]k/n,(k+1)/n]$ respectively solving equation (\ref{Yn}).  
We suppose that there exists $\delta >0$ such that $\sup_n \mathbb{E}[|\xi^n|^{1+\delta}]^{\frac{1}{1+\delta}} < + \infty$ and 
$\sup_n \mathbb{E}[|\tau^n|^{1+\delta}]^{\frac{1}{1+\delta}} < + \infty$ and that we have the 
convergences $\xi^n \xrightarrow{\mathbb{P}} \xi$, $\tau^n \xrightarrow{\mathbb{P}} \tau$ and $W^n \xrightarrow{\mathbb{P}} W$. Then 
$$\forall L \in \mathbb{N}, \ \sup_{t \in [0,L]} |Y^n_{t \wedge \tau^n} - Y_{t \wedge \tau}| 
	+ \int_0^{\tau \wedge \tau^n} |Z^n_{t \wedge \tau^n} - Z_{t \wedge \tau}|^2dt \xrightarrow{\mathbb{P}} 0$$
which we shall denote by $(Y^n, Z^n) \to (Y,Z)$ and also
$$\forall L \in \mathbb{N}, \ \sup_{t \in [0,L]} \left| \int_0^{t \wedge \tau^n} Z^n_s dW^n_s - \int_0^{t \wedge \tau} Z_s dW_s \right| \xrightarrow{\mathbb{P}} 0.$$
\end{theo}

According to this Theorem, we have the convergence in probability of the solutions under the rather strong assumption that the scaled 
random walks converge in probability to the Brownian motion and the random terminal times also converge. Actually, with Donsker's 
Theorem, we have the convergence in law of the scaled random walks to the Brownian motion. If we only have this convergence in law, we 
obtain the following corollary~:

\begin{cor}
\label{hb_cvl}
Let $(Y,Z)$ be the solution of the BSDE (\ref{Y}) and $(Y^n,Z^n)$ be the processes constant on the 
intervals $[k/n,(k+1)/n[$ and $]k/n,(k+1)/n]$ respectively solution of the equation (\ref{Yn}).  
We suppose that $\forall n, \forall k, \varepsilon^n_k=\varepsilon_k$, $\xi = g(W)$ and $\xi^n=g(W^n)$ with $g$ bounded continuous. 
We assume that there exists $\delta >0$ such that $\sup_n \mathbb{E}[|\tau^n|^{1+\delta}]^{\frac{1}{1+\delta}} < + \infty$.
We also suppose that we have the convergence $(W^n,\tau^n) \xrightarrow{\mathcal{L}} (W,\tau)$. 
Then $\left(Y^n_{. \wedge \tau^n}, \int_0^{. \wedge \tau^n} Z^n_s dW^n_s \right)_n$ converges in law to 
$\left(Y_{. \wedge \tau},\int_0^{. \wedge \tau} Z_s dW_s \right)$ for the Skorokhod topology. 
\end{cor}

\demo
According to the Skorokhod representation Theorem, in a space $(\tilde{\Omega},\tilde{\mathcal{A}},\tilde{\mathbb{P}})$, we can find 
$(\tilde{W}^n, \tilde{\tau}^n)$ with the same law as $(W^n,\tau^n)$ and
$(\tilde{W},\tilde{\tau})$ with the same law as $(W,\tau)$ such that 
$(\tilde{W}^n,\tilde{\tau}^n) \xrightarrow{a.s.} (\tilde{W},\tilde{\tau})$.
We denote by $\tilde{\mathcal{F}}$ the natural filtration of $\tilde{W}$ and by $\tilde{\mathcal{F}}^n$ the natural filtrations of the $\tilde{W}^n$. \\
Let us show that $\tilde{\tau}$ is an $\tilde{\mathcal{F}}$-stopping time. 
Let us fix $t \geqslant 0$ and note that $\{\tau \leqslant t\} \in \mathcal{F}_t$ since $\tau$ is a $\mathcal{F}$-stopping time. 
Then, for every $\varepsilon >0$, we can find $h:\mathbb{R}^k \to \mathbb{R}$ measurable and $s_1, \ldots, s_k$ in $[0,t]$ such that
$$\int|\mathbf{1}_{\{\tau \leqslant t\}} - h(W_{s_1}, \ldots, W_{s_k})|d\mathbb{P} < \varepsilon.$$
$(W, \tau) \sim (\tilde{W},\tilde{\tau})$ so 
$\mathbf{1}_{\{\tau \leqslant t\}} - h(W_{s_1}, \ldots, W_{s_k}) \sim \mathbf{1}_{\{\tilde{\tau} \leqslant t\}} - h(\tilde{W}_{s_1},\ldots,\tilde{W}_{s_k})$.
Then, $$\int |\mathbf{1}_{\{\tilde{\tau} \leqslant t\}} - h(\tilde{W}_{s_1}, \ldots, \tilde{W}_{s_k})| d\tilde{\mathbb{P}} 
		= \int|\mathbf{1}_{\{\tau \leqslant t\}} - h(W_{s_1}, \ldots, W_{s_k})|d\mathbb{P}.$$ 
So, for every $\varepsilon >0$, we can find $h:\mathbb{R}^k \to \mathbb{R}$ measurable and $s_1, \ldots, s_k$ in $[0,t]$ such that
$$\int |\mathbf{1}_{\{\tilde{\tau} \leqslant t\}} - h(\tilde{W}_{s_1}, \ldots, \tilde{W}_{s_k})| d\tilde{\mathbb{P}} < \varepsilon.$$ 
Hence, $\{\tilde{\tau} \leqslant t\} \in \tilde{\mathcal{F}}_t$. Then $\tilde{\tau}$ is an $\tilde{\mathcal{F}}$-stopping time. 
On the same way, the $\tilde{\tau}^n$'s are $\tilde{\mathcal{F}}^n$-stopping times.\\
Let $(Y^{'n},Z^{'n})$ be the solution of
$$Y^{'n}_{t}= g(\tilde{W}^n) 
+ \int_{t \wedge \tilde{\tau}^n}^{\tilde{\tau}^n} f(Y^{'n}_{(s \wedge \tilde{\tau}^n)-},Z^{'n}_{s \wedge \tilde{\tau}^n})dA^n_s 
- \int_{t \wedge \tilde{\tau}^n}^{\tilde{\tau}^n} Z^{'n}_{s \wedge \tilde{\tau}^n} d\tilde{W}^n_s$$ 
and $(Y',Z')$ be the solution of
$$Y'_{t \wedge \tilde{\tau}} = g(\tilde{W}) + \int_{t \wedge \tilde{\tau}}^{\tilde{\tau}} f(Y'_s,Z'_s)ds 
	- \int_{t \wedge \tilde{\tau}}^{\tilde{\tau}} Z'_s d\tilde{W}_s, \ t \geqslant 0. $$
All the assumptions of Theorem \ref{theo_Donsker} are filled. So we have the convergence in probability of 
$\left(Y'^n_{. \wedge \tilde{\tau}^n}, \int_0^{. \wedge \tilde{\tau}^n} Z'^n_s d\tilde{W}^n_s \right)_n$ to 
$\left(Y'_{. \wedge \tilde{\tau}}, \int_0^{. \wedge \tilde{\tau}} Z'_s d\tilde{W}_s \right)$. 
Then, denoting by $(Y^n_{. \wedge \tau^n},Z^n_{. \wedge \tau^n})$ the solution of (\ref{Y}) and by $(Y_{. \wedge \tau},Z_{. \wedge \tau})$ 
the solution of (\ref{Yn}), since $(Y^n_{. \wedge \tau^n},Z^n_{. \wedge \tau^n},W^n,\tau^n)$ has the same law as  
$(Y^{'n}_{. \wedge \tilde{\tau}^n},Z^{'n}_{. \wedge \tilde{\tau}^n},\tilde{W}^n,\tilde{\tau}^n)$ and 
$(Y_{. \wedge \tau},Z_{. \wedge \tau},W,\tau)$ has the same law as 
$(Y'_{. \wedge \tilde{\tau}},Z'_{. \wedge \tilde{\tau}},\tilde{W},\tilde{\tau})$, we have the convergence of 
$\left(Y^n_{. \wedge \tau^n}, \int_0^{. \wedge \tau^n} Z^n_s dW^n_s \right)_n$ to 
$\left(Y_{. \wedge \tau},\int_0^{. \wedge \tau} Z_s dW_s \right)$ in distribution for the Skorokhod topology.
\findem

\medskip

We point out that in this corollary, it is necessary to have the joint convergence of $((W^n,\tau^n))_n$ to $(W,\tau)$. 
In section \ref{Donsker_ex}, we shall see an example where the stopping times are hitting times. \\

To prove the first convergence of Theorem \ref{theo_Donsker}, we are going to use the Picard approximations $(Y^p, Z^p)$ and 
$((Y^{n,p}, Z^{n,p}))_n$ defined on the following way :
\begin{eqnarray}
\label{Yp}
Y^{p+1}_{t \wedge \tau} & = & \xi + \int_{t \wedge \tau}^{\tau} f(Y^p_s, Z^p_s) ds 
	- \int_{t \wedge \tau}^{\tau} Z^{p+1}_s dW_s,\ \forall t \geqslant 0,  \\
\label{Ynp}
Y^{n,p+1}_{t} & = & \xi^n + \int_{t \wedge \tau^n}^{\tau^n} f(Y^{n,p}_{s-}, Z^{n,p}_s) dA^n_s 
	- \int_{t \wedge \tau^n}^{\tau^n} Z^{n,p+1}_s dW^n_s, 
\end{eqnarray}
with $Y^0=Z^0=Y^{n,0}=Z^{n,0}=0$ and $A^n_t=\frac{[nt]}{n}$. \\
\indent
We have existence and uniqueness of adapted pairs $(Y^p,Z^p)$ and $((Y^{n,p}, Z^{n,p}))_n$ by induction on $p$, using respectively Theorem 
2.1 of Royer in \cite{article_Manuela} and Proposition \ref{exist_unic_Yn*}. \\

We write :
\begin{eqnarray*}
Y^n-Y=(Y^n-Y^{n,p})+(Y^{n,p}-Y^p)+(Y^p-Y), \\
Z^n-Z=(Z^n-Z^{n,p})+(Z^{n,p}-Z^p)+(Z^p-Z).
\end{eqnarray*}

We are going to prove successively the convergence of the Picard approximations to the solutions, $i.e.$ for every $n$, 
$(Y^{n,p}, Z^{n,p}) \to (Y^n, Z^n)$ and $(Y^p, Z^p) \to (Y, Z)$, and then check that the first convergence of the theorem is true for the 
Picard approximations, $i.e.$ for every $p$, $(Y^{n,p},Z^{n,p}) \to (Y^p,Z^p)$. \\

Most of the proof uses the same arguments as in the proof of Theorem 2.1 of Briand, Delyon and M\'emin in \cite{BDM_Donsker} but there are 
some technical difficulties due to the stopping times. \\

The arguments of Lemma 4.1 in \cite{BDM_Donsker} are still true when we replace a deterministic terminal time by a bounded random terminal
time. So we have :
\begin{equation}
\label{recPic1}
\forall L, \ \sup_n \mathbb{E}\left[\sup_{t \in [0,L]} |Y^n_{\tau^n \wedge t} - Y^{n,p}_{\tau^n \wedge t} |^2 
	+ \int_0^{+\infty} |Z^n_{t} - Z^{n,p}_{t} |^2 dt \right] \xrightarrow[p \to +\infty]{} 0.
\end{equation}

With a truncation argument (using the fact that $\tau$ is almost surely finite), we deduce the convergence of $(Y^p,Z^p)$ to $(Y,Z)$, ie
\begin{equation}
\label{recPic2}
\forall L, \ \mathbb{E}\left[\sup_{t \in [0,L]} |Y^p_{\tau \wedge t} - Y_{\tau \wedge t} |^2 
	+ \int_0^{+\infty} |Z^p_{t} - Z_{t} |^2 dt \right] \xrightarrow[p \to +\infty]{} 0.
\end{equation}

\medskip

Now, let us show that for every $p$, $(Y^{n,p},Z^{n,p}) \to (Y^p,Z^p)$ as $n \to +\infty$ that is 
\begin{equation}
\label{HR}
\forall L \in \mathbb{N}, \ 
\sup_{t \in [0,L]} |Y^{n,p}_{t \wedge \tau^n} - Y^p_{t \wedge \tau}| + 
	\int_0^{\tau^n} |Z^{n,p}_{t \wedge \tau^n} - Z^p_{t \wedge \tau}|^2 dt \xrightarrow{\mathbb{P}} 0 \text{~~as~} n \to \infty.
\end{equation}

We argue by induction on $p$. \\
- The property is true for $p=0$ because $Y^0=Z^0=Y^{n,0}=Z^{n,0}=0$. \\
- We suppose that, for given $p$, (\ref{HR}) holds. Let us prove that (\ref{HR}) is still true for $p+1$. The proof will be given through 
three steps. \\
\indent
In a first step, we introduce the sequence $(M^n)_n$ of processes defined by 
$$M^n_t = Y^{n,p+1}_{t \wedge \tau^n} + \int_0^{t \wedge \tau^n} f(Y^{n,p}_s, Z^{n,p}_s)dA^n_s$$
and the process $M$ defined by 
$$M_t = Y^{p+1}_{t \wedge \tau} + \int_0^{t \wedge \tau} f(Y^{p}_s, Z^{p}_s)ds.$$
In Lemma \ref{hb_lem14}, we prove that $(M^n_{. \wedge \tau^n})_{n}$ is a sequence of 
$(\mathcal{F}^{n, \tau^n})$-martingales. Then, with Lemmas \ref{cvMn} and \ref{lemme_cvM}, we show that we have the following convergence :
$$\forall L, \ \sup_{t \in [0,L]} |M^n_{t \wedge \tau^n} - M_{t \wedge \tau}| \xrightarrow{\mathbb{P}} 0 \text{~~as~} n \to \infty.$$
\indent
The aim of the second step is to prove Lemma \ref{cvZ} :
$$\int_0^{\tau^n} |Z^{n,p+1}_{t \wedge \tau^n} - Z^{p+1}_{t \wedge \tau}|^2 dt \xrightarrow{\mathbb{P}} 0 \text{~~as~} n \to \infty.$$
\indent
At last, in a third step (Lemma \ref{cvY}), we deal with the convergence of $Y^{n,p+1}$ to $Y^{p+1}$ :
$$\forall L, \ \sup_{t \in [0,L]} |Y^{n,p+1}_{t \wedge \tau^n} - Y^{p+1}_{t \wedge \tau}| \xrightarrow{\mathbb{P}} 0 \text{~~as~} n \to \infty.$$

{\bf Step 1 :} Let $(M^n)_n$ be the processes defined by 
$M^n_t = Y^{n,p+1}_{t \wedge \tau^n} + \int_0^{t \wedge \tau^n} f(Y^{n,p}_s, Z^{n,p}_s)dA^n_s$.

\begin{lemme}
\label{hb_lem14}
For each $n$, $M^n$ is a $\mathcal{F}^{n, \tau^n}$-martingale.
\end{lemme}

\demo
Let us show that $M^n_{t} = M^n_0 + \int_0^{t \wedge \tau^n} Z^{n,p+1}_{s \wedge \tau^n}dW^n_s$. 
Using (\ref{Ynp}),
\begin{eqnarray*}
\lefteqn{M^n_0 + \int_0^{t \wedge \tau^n} Z^{n,p+1}_{s \wedge \tau^n}dW^n_s}\\
& = & Y^{n,p+1}_0 + \int_0^{t \wedge \tau^n} Z^{n,p+1}_{s \wedge \tau^n}dW^n_s \\
& = & \xi^n + \int_{0}^{\tau^n} f(Y^{n,p}_{s-}, Z^{n,p}_s) dA^n_s - \int_{0}^{\tau^n} Z^{n,p+1}_s dW^n_s 
	+ \int_0^{t \wedge \tau^n} Z^{n,p+1}_{s \wedge \tau^n}dW^n_s \\
& = & Y^{n,p+1}_{t \wedge \tau^n} + \int_0^{t \wedge \tau^n} f(Y^{n,p}_{s-}, Z^{n,p}_s)dA^n_s = M^n_t .
\end{eqnarray*}

Then, for every $n$, as the process $(M^n_0 + \int_0^{t} Z^{n,p+1}_{s \wedge \tau^n}dW^n_s)_{t \geqslant 0}$ is a $\mathcal{F}^n$-martingale,
the stopped process $(M^n_0 + \int_0^{t \wedge \tau^n} Z^{n,p+1}_{s \wedge \tau^n}dW^n_s)_{t \geqslant 0} = 
(M^n_t)_{t \geqslant 0}$ is a $\mathcal{F}^{n, \tau^n}$-martingale. 
\findem

So, we have $M^n_{t}=\mathbb{E}[M^n_{\tau^n}|\mathcal{F}^{n,\tau^n}_t]$ with
$M^n_{\tau^n}= Y^{n,p+1}_{\tau^n} + \int_0^{\tau^n} f(Y^{n,p}_{s-}, Z^{n,p}_s)dA^n_s$, for every $n$, for every $t$. \\

Before proving the convergence of $(M^n_{\tau^n})_n$, let us see a lemma that links stopped integrals and integrals of stopped 
processes :
\begin{lemme}
\label{int}
For every t, we have the following relations :
\begin{displaymath}
\begin{array}{c}
\int_0^{t \wedge \tau} f(Y^{p}_s,Z^{p}_s)ds = \int_0^t f(Y^{p}_{s \wedge \tau},Z^{p}_{s \wedge \tau})ds + 
	(\tau \wedge t - t)f(Y^{p}_\tau,Z^{p}_\tau), \\
\int_0^{t \wedge \tau^n} f(Y^{n,p}_s,Z^{n,p}_s)ds = \int_0^t f(Y^{n,p}_{s \wedge \tau^n},Z^{n,p}_{s \wedge \tau^n})ds + 
	(\tau^n \wedge t - t)f(Y^{n,p}_{\tau^n},Z^{n,p}_{\tau^n}).
\end{array}
\end{displaymath} 
\end{lemme}
\demo
Let us prove the first relation. \\
On $\{t \leqslant \tau\}$,
$$\int_0^t f(Y^{p}_{s \wedge \tau},Z^{p}_{s \wedge \tau})ds
	= \int_0^t f(Y^{p}_{s},Z^{p}_{s})ds
	= \int_0^{t \wedge \tau} f(Y^{p}_{s},Z^{p}_{s})ds$$
and on $\{t > \tau\}$, we have :
\begin{eqnarray*}
\int_0^t f(Y^{p}_{s \wedge \tau},Z^{p}_{s \wedge \tau})ds 
& = & \int_0^{t \wedge \tau} f(Y^{p}_{s \wedge \tau},Z^{p}_{s \wedge \tau})ds 
	+ \int_{t \wedge \tau}^t f(Y^{p}_{s \wedge \tau},Z^{p}_{s \wedge \tau})ds \\
& = & \int_0^{t \wedge \tau} f(Y^{p}_{s},Z^{p}_{s})ds + \int_{t \wedge \tau}^t f(Y^{p}_{\tau},Z^{p}_{\tau})ds \\
& = & \int_0^{t \wedge \tau} f(Y^{p}_{s},Z^{p}_{s})ds + (t - t \wedge \tau) f(Y^{p}_{\tau},Z^{p}_{\tau}). 
\end{eqnarray*}
The second equality is proved by the same arguments.
\findem

\begin{lemme}
\label{cvMn}
$(M^n_{\tau^n})_n$ converges in $L^1$ to $Y^{p+1}_{\tau} + \int_0^{\tau} f(Y^{p}_s, Z^{p}_s)ds$.
\end{lemme}

\demo
\indent
It suffices to show that $(M^n_{\tau^n})_n$ converges to $Y^{p+1}_{\tau} + \int_0^{\tau} f(Y^{p}_s, Z^{p}_s)ds$ in probability. 
Indeed, for every $n$, 
$\mathbb{E}[|M^n_{\tau^n}|^{1+\delta}]^{\frac{1}{1+\delta}} 
	\leqslant \mathbb{E}[|\xi^n|^{1+\delta}]^{\frac{1}{1+\delta}} +  \|f\|_{\infty}\mathbb{E}[|\tau^n|^{1+\delta}]^{\frac{1}{1+\delta}}$. 
So, when we take the $\sup$ in $n$, $\sup_n \mathbb{E}[|M^n_{\tau^n}|^{1+\delta}]^{\frac{1}{1+\delta}} < \infty$
according to the assumptions on $(\xi^n)_n$, $(\tau^n)_n$ and $f$. 
So, the sequence $(M^n_{\tau^n})$ is uniformly integrable and then the convergence in probability implies the convergence in $L^1$.\\

Let us show that $(M^n_{\tau^n})_n$ converges in probability to $Y^{p+1}_{\tau} + \int_0^{\tau} f(Y^{p}_s, Z^{p}_s)ds$.
\begin{eqnarray*}
\lefteqn{\left|M^n_{\tau^n} - Y^{p+1}_{\tau} + \int_0^{\tau} f(Y^{p}_s, Z^{p}_s)ds \right|}\\
& \leqslant & |Y^{n,p+1}_{\tau^n} - Y^{p+1}_{\tau}|
	+ \left|\int_0^{\tau^n} f(Y^{n,p}_{s-}, Z^{n,p}_s)dA^n_s - \int_0^{\tau} f(Y^{p}_s, Z^{p}_s)ds \right|
\end{eqnarray*}
But, 
\begin{equation}
\label{cv1}
|Y^{n,p+1}_{\tau^n} - Y^{p+1}_{\tau}| = |\xi^n - \xi| \xrightarrow{\mathbb{P}} 0.
\end{equation}

We are going to conclude by showing that 
$$\left|\int_0^{\tau^n} f(Y^{n,p}_{s-}, Z^{n,p}_s)dA^n_s - \int_0^{\tau} f(Y^{p}_s, Z^{p}_s)ds \right| \xrightarrow{\mathbb{P}} 0.$$

For each $n$, for each $\omega$, there exists a unique $k_{n}$ such that $\tau^{n}$ is in the interval $[k_n/n,(k_{n}+1)/n[$. As the 
processes $Y^{n,p}$ and $Z^{n,p}$ are constant on the intervals of the form $[k/n,(k+1)/n[$ and $]k/n,(k+1)/n]$ respectively by 
construction, we have :
$$\int_0^{\tau^{n}} f(Y^{n,p}_{s-}, Z^{n,p}_s)dA^{n}_s 	= \int_0^{k_{n}/n} f(Y^{n,p}_s, Z^{n,p}_s)ds.$$
So, 
\begin{eqnarray*}
\lefteqn{\left|\int_0^{\tau^{n}} f(Y^{n,p}_{s-}, Z^{n,p}_s)dA^{n}_s - \int_0^{\tau} f(Y^{p}_s, Z^{p}_s)ds \right|}\\
& = & \left|\int_0^{k_{n}/n} f(Y^{n,p}_s, Z^{n,p}_s)ds - \int_0^{\tau} f(Y^{p}_s, Z^{p}_s)ds\right| \\
& \leqslant & 
\left|\int_0^{(k_{n}/n) \wedge \tau} (f(Y^{n,p}_s, Z^{n,p}_s)-f(Y^{p}_s, Z^{p}_s))ds\right| \\
&& + \Bigg|\int_{(k_{n}/n) \wedge \tau}^{(k_{n}/n) \vee \tau} (f(Y^{n,p}_s, Z^{n,p}_s)\mathbf{1}_{\{\tau \leqslant k_{n}/n\}}  
	-f(Y^{p}_s, Z^{p}_s)\mathbf{1}_{\{\tau > k_{n}/n\}})ds\Bigg|.
\end{eqnarray*}
Taking $t=(k_{n}/n) \wedge \tau$ in the two equalities of Lemma \ref{int}, we have :
\begin{eqnarray*}
\int_0^{(k_{n}/n) \wedge \tau} f(Y^{p}_s, Z^{p}_s) ds 
	& = & \int_0^{(k_{n}/n) \wedge \tau} f(Y^{p}_{s \wedge \tau}, Z^{p}_{s \wedge \tau}) ds,\\
\int_0^{(k_{n}/n) \wedge \tau} f(Y^{n,p}_s, Z^{n,p}_s) ds
	& = & \int_0^{(k_{n}/n) \wedge \tau} f(Y^{n,p}_{s \wedge \tau^{n}}, Z^{n,p}_{s \wedge \tau^{n}}) ds.
\end{eqnarray*}
So,
\begin{eqnarray*}
\lefteqn{\left|\int_0^{(k_{n}/n) \wedge \tau} (f(Y^{n,p}_s, Z^{n,p}_s)-f(Y^{p}_s, Z^{p}_s))ds\right|}\\
& = & \left|\int_0^{(k_{n}/n) \wedge \tau} 
	(f(Y^{n,p}_{s \wedge \tau^{n}}, Z^{n,p}_{s \wedge \tau^{n}}) -f(Y^{p}_{s \wedge \tau}, Z^{p}_{s \wedge \tau}))ds  \right| \\
& \leqslant & \int_0^{(k_{n}/n) \wedge \tau} 
	K[|Y^{n,p}_{s \wedge \tau^{n}} - Y^{p}_{s \wedge \tau}| + |Z^{n,p}_{s \wedge \tau^{n}} - Z^{p}_{s \wedge \tau}|]ds 
	\text{~~because $f$ is $K$-Lipschitz} \\
& \leqslant & K \int_0^{\tau^n \wedge \tau} |Y^{n,p}_{s \wedge \tau^{n}} - Y^{p}_{s \wedge \tau}| ds
	+ K \int_0^{\tau \wedge \tau^{n}} |Z^{n,p}_{s \wedge \tau^{n}} - Z^{p}_{s \wedge \tau}|ds.  
\end{eqnarray*}
Using the induction assumption (\ref{HR}), the fact that $\tau < + \infty$ $a.s.$ and Cauchy-Schwarz inequality, we prove that 
\begin{equation}
\label{cv_intZ}
\int_0^{\tau \wedge \tau^{n}} |Z^{n,p}_{s \wedge \tau^{n}} - Z^{p}_{s \wedge \tau}|^2ds \xrightarrow{\mathbb{P}} 0.
\end{equation}
On the other hand, let us fix $\varepsilon > 0$ and $\eta > 0$. Since $\tau < + \infty$ $a.s.$, we can find $T$ such that 
$\mathbb{P}[\tau \geqslant T] \leqslant \varepsilon.$ Then,
\begin{eqnarray*}
\lefteqn{\mathbb{P}\left[\int_0^{\tau^n \wedge \tau} |Y^{n,p}_{s \wedge \tau^n} - Y^p_{s \wedge \tau}| ds \geqslant \eta \right]}\\
& = & \mathbb{P}\left[\int_0^{\tau^n \wedge \tau} |Y^{n,p}_{s \wedge \tau^n} - Y^p_{s \wedge \tau}| ds \ \mathbf{1}_{\{\tau \wedge \tau^n < T\}} 
		\geqslant \eta \right] 
+ \mathbb{P}\left[\int_0^{\tau^n \wedge \tau} |Y^{n,p}_{s \wedge \tau^n} - Y^p_{s \wedge \tau}| ds \ \mathbf{1}_{\{\tau \wedge \tau^n \geqslant T\}} 
		\geqslant \eta \right] \\
& \leqslant & \mathbb{P} \left[ T \sup_{s \in [0,T]} |Y^{n,p}_{s \wedge \tau^n} - Y^p_{s \wedge \tau}| \geqslant \eta \right] 
+ \mathbb{P}[\tau \wedge \tau^n \geqslant T] \\
& \leqslant & 2\varepsilon \text{~~by choice of $T$ and using (\ref{HR}).}
\end{eqnarray*}
So,
\begin{equation}
\label{cv_intY}
\int_0^{\tau^n \wedge \tau} |Y^{n,p}_{s \wedge \tau^n} - Y^p_{s \wedge \tau}| ds \xrightarrow{\mathbb{P}} 0.
\end{equation}
Finally, using the convergences (\ref{cv_intZ}) and (\ref{cv_intY}), we have
\begin{equation}
\label{A}
\left|\int_0^{(k_{n}/n) \wedge \tau} (f(Y^{n,p}_s, Z^{n,p}_s)-f(Y^{p}_s, Z^{p}_s))ds\right| \xrightarrow{\mathbb{P}} 0.
\end{equation}
On the other hand, $|\tau^{n} - k_{n}/n| \leqslant 1/n$ and $\tau^{n} \xrightarrow{\mathbb{P}} \tau$, so we have $(k_{n}/n) \xrightarrow{\mathbb{P}} \tau$. Then, 
\begin{eqnarray}
\label{B}
\lefteqn{\left|\int_{(k_{n}/n) \wedge \tau}^{(k_{n}/n) \vee \tau} 
  (f(Y^{n,p}_s, Z^{n,p}_s){\bf 1}_{\{\tau \leqslant k_{n}/n\}}-f(Y^{p}_s, Z^{p}_s){\bf 1}_{\{\tau > k_{n}/n\}}) ds \right| }
   \nonumber \\
& \leqslant & |k_{n}/n - \tau| \ \|f\|_{\infty}
\quad \quad \quad \quad  \quad \quad \quad \quad  \quad \quad \quad \quad \quad \quad  \quad \quad \quad \quad \quad  \quad \quad\\
&  \xrightarrow{\mathbb{P}} & 0. \nonumber
\end{eqnarray}
According to the convergences (\ref{A}) and (\ref{B}), we have 
\begin{equation}
\label{cvf}
\left|\int_0^{\tau^{n}} f(Y^{n,p}_{s-}, Z^{n,p}_s)dA^{n}_s - \int_0^{\tau} f(Y^{p}_s, Z^{p}_s)ds\right| \xrightarrow{\mathbb{P}} 0.
\end{equation}
Hence, according to the convergences (\ref{cv1}) and (\ref{cvf}),
\begin{equation}
\label{hb_cv_Mn}
M^n_{\tau^n} \xrightarrow{\mathbb{P}} \left(Y^{p+1}_{\tau} + \int_0^{\tau} f(Y^{p}_s, Z^{p}_s)ds\right) .
\end{equation}
At last, Lemma \ref{cvMn} is proved. 
\findem

Let $M$ be the process defined by $M_t=Y^{p+1}_{t \wedge \tau} + \int_0^{t \wedge \tau} f(Y^{p}_s, Z^{p}_s)ds$. 

\begin{lemme}
\label{lemme_cvM}
We have the convergence 
\begin{equation}
\label{cvM}
\forall L, \ \sup_{t \in [0,L]} |M^n_{t} - M_{t}| \xrightarrow{\mathbb{P}} 0.
\end{equation}
\end{lemme}

\demo
Using the same computation as for $M_{t}^n$ in Lemma \ref{hb_lem14}, we find :
$$M_{t} = M_0 + \int_0^{t \wedge \tau} Z_s^{p+1} dW_s = Y^{p+1}_{t \wedge \tau} + \int_0^{t \wedge \tau} f(Y^{p}_s, Z^{p}_s)ds.$$
So the process $M$ is a $\mathcal{F}^\tau$-martingale. 
Then, $M_{t \wedge \tau} = \mathbb{E}[M_\tau|\mathcal{F}^\tau_t]$, for every $t$, where $M_\tau = Y^{p+1}_{\tau} + \int_0^{\tau} f(Y^{p}_s, Z^{p}_s)ds$. 
According to (\ref{hb_cv_Mn}), $M^n_{\tau^n} \xrightarrow{L^1} M_\tau$. \\

If we prove that we have the convergence of the sequence of filtrations $(\mathcal{F}^{n,\tau^n})_n$ to $\mathcal{F}^\tau$, according to Remark 1.2 in 
Coquet, M\'emin and S\l omi\'nski \cite{cvfiltration}, we have $M^n_{.} \xrightarrow{\mathbb{P}} M_{.}$ for the Skorokhod topology.\\
\indent 
$(W^n)_n$ is a sequence of processes with independant increments that converges in probability to $W$. So, according to Proposition 2 in 
\cite{cvfiltration}, $\mathcal{F}^n \xrightarrow{w} \mathcal{F}$. Moreover, $W$ is continuous and $\tau^n \xrightarrow{\mathbb{P}} \tau$ so
according to Corollary \ref{cv_proc_arrete_continu}, $\mathcal{F}^{n,\tau^n} \xrightarrow{w} \mathcal{F}^\tau$. \\

Finally, we have $M^n_{.} \xrightarrow{\mathbb{P}} M_{.}$ for the Skorokhod topology. Moreover, $\mathcal{F}$ is a Brownian filtration, so, 
every $\mathcal{F}$-martingale is continuous. In particular, $\mathbb{E}[M_\tau|\mathcal{F}_.]$ is continuous. The stopped process is also continuous, 
$i.e.$ $M_{.} = \mathbb{E}[M_\tau|\mathcal{F}^\tau_.]$ is a continuous process. So, the previous convergence is uniform in $t$ on every compact set and 
Lemma \ref{lemme_cvM} is proved.
\findem

\medskip

{\bf Step 2 :} In this step, we shall prove the convergence of $(Z^{n,p+1})$ to $Z^{p+1}$ as $n \to \infty$. \\

First, let us see a lemma of convergence of quadratic variations :
\begin{lemme}
\label{cv_crochet}
$(W^{n,\tau^n}, M^{n,\tau^n}, [M^{n,\tau^n},M^{n,\tau^n}], [M^{n,\tau^n},W^{n,\tau^n}])$ converges to  
$(W^\tau, M^\tau, [M^\tau,M^\tau], [M^\tau,W^\tau])$ in probability.
\end{lemme}
\demo
It suffices to follow the lines of Briand, Delyon and M\'emin in the proof of Theorem 3.1 in \cite{BDM_Donsker}. 
\findem

\begin{lemme}
\label{cvZ}
$\int_0^{\tau^n} |Z^{n,p+1}_{t \wedge \tau^n} - Z^{p+1}_{t \wedge \tau}|^2 ds \xrightarrow{\mathbb{P}} 0$.
\end{lemme}
\demo
To prove this lemma, we prove that we have these two convergences :
\begin{eqnarray}
\label{hb_cvL2}
\int_0^{\tau^n} |Z^{n,p+1}_{t \wedge \tau^n}|^2 dt - \int_0^{\tau^n} |Z^{p+1}_{t \wedge \tau}|^2 dt \xrightarrow{\mathbb{P}} 0,\\
\label{hb_cvL2f}
\forall g \in L^2, \int_0^{\tau^n} g(s)(Z^{n,p+1}_{s \wedge \tau^n} - Z^{p+1}_{s \wedge \tau})ds \xrightarrow{\mathbb{P}} 0.
\end{eqnarray}
\noindent Then, we conclude with the following lines :
\begin{eqnarray*}
\lefteqn{\int_0^{\tau^n} |Z^{n,p+1}_{t \wedge \tau^n} - Z^{p+1}_{t \wedge \tau}|^2 dt} \\
& = & \int_0^{\tau^n} |Z^{n,p+1}_{t \wedge \tau^n}|^2 dt + \int_0^{\tau^n} |Z^{p+1}_{t \wedge \tau}|^2 dt
-2 \int_0^{\tau^n} Z^{n,p+1}_{t \wedge \tau^n}Z^{p+1}_{t \wedge \tau} dt \\
& = & \left(\int_0^{\tau^n} |Z^{n,p+1}_{t \wedge \tau^n}|^2 dt -  \int_0^{\tau^n} |Z^{p+1}_{t \wedge \tau}|^2 dt \right) 
 + 2 \left( \int_0^{\tau^n} |Z^{p+1}_{t \wedge \tau}|^2 dt - \int_0^{\tau^n} Z^{n,p+1}_{t \wedge \tau^n}Z^{p+1}_{t \wedge \tau} dt \right) \\
& \xrightarrow{\mathbb{P}} & 0 \text{~~according to (\ref{hb_cvL2}) and (\ref{hb_cvL2f}) with $g=Z^{p+1}$}.
\end{eqnarray*}

First, let us show that 
$\int_0^{\tau^n} |Z^{n,p+1}_{t \wedge \tau^n}|^2 dt - \int_0^{\tau^n} |Z^{p+1}_{t \wedge \tau}|^2 dt \xrightarrow{\mathbb{P}} 0$. \\
\begin{eqnarray*}
\lefteqn{\left| \int_0^{\tau^n} |Z^{n,p+1}_{t \wedge \tau^n}|^2 dt - \int_0^{\tau^n} |Z^{p+1}_{t \wedge \tau}|^2 dt \right|} \\
& \leqslant & \left| \int_0^{\tau^n} |Z^{n,p+1}_{t \wedge \tau^n}|^2 dt - \int_0^{\tau^n} |Z^{n,p+1}_{t \wedge \tau^n}|^2 dA^n_t \right|
+ \left| \int_0^{\tau^n} |Z^{n,p+1}_{t \wedge \tau^n}|^2 dA^n_t - \int_0^{\tau^n} |Z^{p+1}_{t \wedge \tau}|^2 dt \right|.
\end{eqnarray*}
According to Lemma \ref{cv_crochet}, 
$$\forall L, \ \sup_{t \in [0,L]} \left|[M^{n,\tau^n},M^{n,\tau^n}]_t - [M^\tau,M^\tau]_t \right| \xrightarrow{\mathbb{P}} 0.$$ 
But,
$$M^{n,\tau^n}_t=M^n_0 + \int_0^{t \wedge \tau^n} Z^{n,p+1}_s dW^n_s = M^n_0 
	+ \frac{1}{\sqrt{n}}\sum_{k=1}^{[n(t \wedge \tau^n)]} Z^{n,p+1}_{k/n} \varepsilon^n_k.$$
So, 
$$[M^{n,\tau^n},M^{n,\tau^n}]_t = \frac{1}{n}\sum_{k=1}^{[n(t \wedge \tau^n)]} |Z^{n,p+1}_{k/n}|^2 
	= \int_0^{t \wedge \tau^n} |Z^{n,p+1}_s|^2 dA^n_s.$$
On the other hand, $M^\tau$ is a continuous martingale, so 
$$[M^\tau,M^\tau]_t = <M^\tau,M^\tau>_t = \int_0^{t \wedge \tau} |Z^{p+1}_s|^2 ds.$$
Then,
$$\forall L, \ 
\sup_{s \in [0,L]} \left| \int_0^{s \wedge \tau^n} |Z^{n,p+1}_t|^2 dA^n_t - \int_0^{s \wedge \tau} |Z^{p+1}_t|^2 dt \right| \xrightarrow{\mathbb{P}} 0.$$
In particular, using the fact that $\sup_n \tau^n$ is almost surely finite because $(\tau^n)_n$ converges in probability to the almost
surely finite stopping time $\tau$, we have
$$\left| \int_0^{\tau^n} |Z^{n,p+1}_t|^2 dA^n_t - \int_0^{\tau \wedge \tau^n} |Z^{p+1}_t|^2 dt \right| \xrightarrow{\mathbb{P}} 0.$$
As $\tau^n \xrightarrow{\mathbb{P}} \tau$, by dominated convergence, we have :
$$\left| \int_0^{\tau^n} |Z_t^{p+1}|^2 dt - \int_0^{\tau^n \wedge \tau} |Z^{p+1}_t|^2dt \right| \xrightarrow{\mathbb{P}} 0 \text{~as~} n \to + \infty.$$
Moreover, $\int_0^{\tau^n} |Z^{n,p+1}_t|^2 dA^n_t = \int_0^{\tau^n} |Z^{n,p+1}_t|^2 dt.$ So,  the convergence (\ref{hb_cvL2}) is proved. \\

It remains to prove convergence (\ref{hb_cvL2f}). Let us fix $g \in L^2(\mathbb{R}^+)$ and show that 
$\int_0^{\tau^n} g(s)(Z^{n,p+1}_{s \wedge \tau^n} - Z^{p+1}_{s \wedge \tau})ds \xrightarrow{\mathbb{P}} 0$. \\
According to Lemma \ref{cv_crochet}, we have the convergence 
$$\forall L, \ \sup_{t \in [0,L]} \left|[M^{n,\tau^n},W^{n,\tau^n}]_t - [M^\tau,W^\tau]_t \right| \xrightarrow{\mathbb{P}} 0.$$
But, 
$$[M^{n,\tau^n},W^{n,\tau^n}]_t = \frac{1}{n}\sum_{k=1}^{[n(t \wedge \tau^n)]} Z^{n,p+1}_{k/n} 
	= \int_0^{t \wedge \tau^n} Z^{n,p+1}_s dA^n_s.$$
On the other hand, $M^\tau$ and $W^\tau$ are continuous martingales, so 
$$[M^\tau,W^\tau]_t = <M^\tau,W^\tau>_t = \int_0^{t \wedge \tau} Z^{p+1}_s ds.$$
Then,
\begin{equation}
\label{cv_sup}
\forall L, \ \sup_{t \in [0,L]} \left| \int_0^{t \wedge \tau^n} Z^{n,p+1}_s dA^n_s - \int_0^{t \wedge \tau} Z^{p+1}_s ds \right| \xrightarrow{\mathbb{P}} 0.
\end{equation}
We write :
\begin{eqnarray*}
\lefteqn{\left| \int_0^{\tau \wedge \tau^{n}} g(s)(Z^{n,p+1}_{s \wedge \tau^{n}} - Z^{p+1}_{s \wedge \tau})ds \right|} \\
& \leqslant & \left| \int_0^{\tau \wedge \tau^{n}} g(s)Z^{n,p+1}_{s \wedge \tau^{n}}ds 
		- \int_0^{\tau \wedge \tau^{n}} g(s)Z^{n,p+1}_{s \wedge \tau^{n}} dA^{n}_s \right| \\
&& \quad \quad \quad 
	+ \left| \int_0^{\tau \wedge \tau^{n}} g(s)Z^{n,p+1}_{s \wedge \tau^{n}} dA^{n}_s 
		- \int_0^{\tau \wedge \tau^{n}} g(s) Z^{p+1}_{s \wedge \tau} ds \right|
\end{eqnarray*}
As the function $g$ is measurable and $\tau$ is almost surely finite, using a density argument and the 
convergence (\ref{cv_sup}), we have :
\begin{equation}
\label{hb_cvc2}
\left| \int_0^{\tau^{n}} g(s)Z^{n,p+1}_{s \wedge \tau^{n}}dA^n_s 
	- \int_0^{\tau \wedge \tau^{n}} g(s)Z^{p+1}_{s \wedge \tau^{n}} ds \right|
		\to 0.
\end{equation}
By dominated convergence, we deduce that 
$$\left| \int_0^{\tau \wedge \tau^{n}} g(s)Z^{n,p+1}_{s \wedge \tau^{n}}ds 
		- \int_0^{\tau \wedge \tau^{n}} g(s)Z^{n,p+1}_{s \wedge \tau^{n}} dA^{n}_s \right| \xrightarrow{\mathbb{P}} 0.$$
Next, using the density of the set of continuous function with compact support in $L^2$, we show that 
$$\left| \int_0^{\tau \wedge \tau^{n}} g(s)Z^{n,p+1}_{s \wedge \tau^{n}} dA^{n}_s 
		- \int_0^{\tau \wedge \tau^{n}} g(s) Z^{p+1}_{s \wedge \tau} ds \right| \xrightarrow{\mathbb{P}} 0.$$
The convergence (\ref{hb_cvL2f}) is now proved.

\noindent Lemma \ref{cvZ} is proved.
\findem

{\bf Step 3 :} Let us prove the convergence of $Y^{n,p+1}$ to $Y^{p+1}$ when $n$ goes to $\infty$.
\begin{lemme}
\label{cvY}
For every $L$, ${\displaystyle \sup_{t \in [0,L]} |Y^{n,p+1}_{t \wedge \tau^n} - Y^{p+1}_{t \wedge \tau}| \xrightarrow{\mathbb{P}} 0}$. 
\end{lemme}

\demo
We fix $L$. We write :
\begin{eqnarray*}
\lefteqn{\sup_{t \in [0,L]} |Y^{n,p+1}_{t \wedge \tau^n} - Y^{p+1}_{t \wedge \tau}|} \\
& \leqslant & \sup_{t \in [0,L]} |M^n_{t \wedge \tau^n} - M_{t \wedge \tau}| 
+ \sup_{t \in [0,L]} \left| \int_0^{t \wedge \tau^n} f(Y^{n,p}_{s-},Z^{n,p}_s)dA^n_s - \int_0^{t \wedge \tau} f(Y^{p}_s,Z^{p}_s)ds \right|
\end{eqnarray*}
Moreover,
\begin{eqnarray*}
\lefteqn{\sup_{t \in [0,L]} \left| \int_0^{t \wedge \tau^n} f(Y^{n,p}_{s-},Z^{n,p}_s)dA^n_s - \int_0^{t \wedge \tau} f(Y^{p}_s,Z^{p}_s)ds \right|} \\
& \leqslant & 
\sup_{t \in [0,L]} \left| \int_0^{t \wedge \tau^n} f(Y^{n,p}_{s-},Z^{n,p}_s)dA^n_s - \int_0^{t \wedge \tau^n} f(Y^{n,p}_s,Z^{n,p}_s)ds \right| \\
&& + 
\sup_{t \in [0,L]} \left| \int_0^{t \wedge \tau^n} f(Y^{n,p}_s,Z^{n,p}_s)ds - \int_0^{t \wedge \tau^n} f(Y^{p}_s,Z^{p}_s)ds \right|
		{\bf 1}_{{\tau^n \leqslant \tau}} \\
&& + 
\sup_{t \in [0,L]} \left| \int_0^{t \wedge \tau} f(Y^{n,p}_s,Z^{n,p}_s)ds - \int_0^{t \wedge \tau} f(Y^{p}_s,Z^{p}_s)ds \right|
		{\bf 1}_{{\tau^n > \tau}} \\
&& + 
\sup_{t \in [0,L]} \left| \int_0^{t \wedge \tau^n} f(Y^{p}_s,Z^{p}_s)ds -  \int_0^{t \wedge \tau} f(Y^{p}_s,Z^{p}_s)ds \right|
		{\bf 1}_{{\tau^n \leqslant \tau}}\\
&& + 
\sup_{t \in [0,L]} \left| \int_0^{t \wedge \tau^n} f(Y^{n,p}_s,Z^{n,p}_s)ds -  \int_0^{t \wedge \tau} f(Y^{n,p}_s,Z^{n,p}_s)ds \right|
		{\bf 1}_{{\tau^n > \tau}}
\end{eqnarray*}
But, for every $t$, for every $\omega$, we can find $k_n$ such that $t \wedge \tau^n\in [k_n/n, (k_n+1)/n[$. So, we have
\begin{eqnarray*}
\left| \int_0^{t \wedge \tau^n} f(Y^{n,p}_{s-},Z^{n,p}_s)dA^n_s - \int_0^{t \wedge \tau^n} f(Y^{n,p}_s,Z^{n,p}_s)ds \right| 
& = & \left| \int_{k_n/n}^{t \wedge \tau^n} f(Y^{n,p}_s,Z^{n,p}_s)ds \right| \\
& \leqslant & |t \wedge \tau^n - k_n/n| \ \|f\|_{\infty} \\
& \leqslant & \|f\|_{\infty}/n.
\end{eqnarray*}
Then,
\begin{equation}
\label{cvA}
\sup_{t \in [0,L]} 
	\left| \int_0^{t \wedge \tau^n} f(Y^{n,p}_{s-},Z^{n,p}_s)dA^n_s - \int_0^{t \wedge \tau^n} f(Y^{n,p}_s,Z^{n,p}_s)ds \right|
\leqslant  \|f\|_{\infty}/n \xrightarrow{\mathbb{P}} 0.
\end{equation}
Next,
\begin{eqnarray}
\label{cvB}
\lefteqn{\sup_{t \in [0,L]} \left| \int_0^{t \wedge \tau^n} f(Y^{n,p}_s,Z^{n,p}_s)ds 
		- \int_0^{t \wedge \tau^n} f(Y^{p}_s,Z^{p}_s)ds \right|{\bf 1}_{\{\tau^n \leqslant \tau\}}} \nonumber \\
&& \quad \quad \quad 	 + \sup_{t \in [0,L]} \left| \int_0^{t \wedge \tau} f(Y^{n,p}_s,Z^{n,p}_s)ds 
		- \int_0^{t \wedge \tau} f(Y^{p}_s,Z^{p}_s)ds \right|{\bf 1}_{\{\tau^n > \tau\} } \nonumber \\
& = & 
\sup_{t \in [0,L]} \left| \int_0^{t \wedge \tau^n} f(Y^{n,p}_{s \wedge \tau^n},Z^{n,p}_{s \wedge \tau^n})ds 
	- \int_0^{t \wedge \tau^n} f(Y^{p}_{s \wedge \tau},Z^{p}_{s \wedge \tau})ds \right|{\bf 1}_{\{\tau^n \leqslant \tau\}} \nonumber \\
&& \quad \quad \quad + \sup_{t \in [0,L]} \left| \int_0^{t \wedge \tau} f(Y^{n,p}_{s \wedge \tau^n},Z^{n,p}_{s \wedge \tau^n})ds 
	- \int_0^{t \wedge \tau} f(Y^{p}_{s \wedge \tau},Z^{p}_{s \wedge \tau})ds \right|{\bf 1}_{\{\tau^n > \tau\}}  \\
& \leqslant &
  \int_0^{\tau \wedge \tau^n \wedge L} 
  	|f(Y^{n,p}_{s \wedge \tau^n},Z^{n,p}_{s \wedge \tau^n}) - f(Y^{p}_{s \wedge \tau},Z^{p}_{s \wedge \tau})|ds
   \nonumber \\
& \leqslant & K(\tau \wedge \tau^n \wedge L)\sup_{t \in [0,L]} |Y^{n,p}_{t \wedge \tau^n} - Y^{p}_{t \wedge \tau}|
	+ K(\tau \wedge \tau^n \wedge L) \left( \int_0^{\tau \wedge \tau^n} |Z^{n,p}_{s \wedge \tau^n} - Z^{p}_{s \wedge \tau}|^2 ds \right)^{1/2}  \nonumber \\
&&  \quad \quad \text{~~because $f$ is $K$-Lipschitz and using Cauchy-Schwarz inequality}   \nonumber  \\
& \xrightarrow{\mathbb{P}} & 0 \text{~~by induction assumption (\ref{HR}).} \nonumber
\end{eqnarray}
Moreover,
\begin{eqnarray}
\label{cvC}
\lefteqn{\sup_{t \in [0,L]} \left| \int_0^{t \wedge \tau^n} f(Y^{p}_s,Z^{p}_s)ds -  \int_0^{t \wedge \tau} f(Y^{p}_s,Z^{p}_s)ds \right|
	{\bf 1}_{\{\tau^n \leqslant \tau\}}} \nonumber \\
&& \quad \quad
+ \sup_{t \in [0,L]} \left| \int_0^{t \wedge \tau^n} f(Y^{n,p}_s,Z^{n,p}_s)ds -  \int_0^{t \wedge \tau} f(Y^{n,p}_s,Z^{n,p}_s)ds \right|
{\bf 1}_{\{\tau^n >\tau\}} \nonumber \\
& \leqslant & \sup_{t \in [0,L]} |t \wedge \tau^n - t \wedge \tau| \ \|f\|_{\infty} \\
& \xrightarrow{\mathbb{P}} & 0 \text{~~because~~} \tau^n \xrightarrow{\mathbb{P}} \tau. \nonumber
\end{eqnarray}

Finally, according to the convergences (\ref{cvM}), (\ref{cvA}), (\ref{cvB}) and (\ref{cvC}), we have the convergence anounced in Lemma 
\ref{cvY}.
\findem

Then, using Lemma \ref{cvZ} and Lemma \ref{cvY}, we have
$$\forall L, \ \sup_{t \in [0,L]} |Y^{n,p+1}_{t \wedge \tau^n} - Y^{p+1}_{t \wedge \tau}| 
	+ \int_0^{\tau^n} |Z^{n,p+1}_{t \wedge \tau^n} - Z^{p+1}_{t \wedge \tau}|^2 dt \xrightarrow{\mathbb{P}} 0.$$
That concludes the proof of the induction. \\

We have shown that (\ref{HR}) is true for every $p$, $i.e.$
\begin{equation}
\label{recurrence}
\forall p, \ \forall L, \ \sup_{t \in [0,L]} |Y^{n,p}_{t \wedge \tau^n} - Y^{p}_{t \wedge \tau}| 
	+ \int_0^{\tau^n} |Z^{n,p}_{t \wedge \tau^n} - Z^{p}_{t \wedge \tau}|^2 dt \xrightarrow{\mathbb{P}} 0.
\end{equation}

Thanks to convergences (\ref{recPic1}), (\ref{recPic2}) and (\ref{recurrence}), the first part of Theorem \ref{theo_Donsker} is proved, $i.e.$ 
$$\forall L, \ \sup_{t \in [0,L]} |Y^{n}_{t \wedge \tau^n} - Y_{t \wedge \tau}| 
	+ \int_0^{\tau^n} |Z^{n}_{t \wedge \tau^n} - Z_{t \wedge \tau}|^2 dt \xrightarrow{\mathbb{P}} 0.$$

To prove the second part of the theorem, we define the processes $M^n$ and $M$ by 
$$M^n_t = M^n_0 + \int_0^{t \wedge \tau^n} Z^n_{s \wedge \tau^n} dW^n_s = Y^n_{t \wedge \tau^n} + \int_0^{t \wedge \tau^n}
f(Y^n_s,Z^n_s)dA^n_s, \forall t$$
and 
$$M_t = M_0 + \int_0^{t \wedge \tau} Z_{s \wedge \tau} dW_s = Y_{t \wedge \tau} + \int_0^{t \wedge \tau} f(Y_s,Z_s)ds, \forall t.$$
It is clear that, for each $n$, $M^n$ is a $\mathcal{F}^{n,\tau^n}$-martingale and that $M$ is a $\mathcal{F}^\tau$-martingale. So, we can write
$$M^n_t=\mathbb{E}[M_{\tau^n} | \mathcal{F}^{n,\tau^n}_t] \text{~and~} M_t=\mathbb{E}[M_\tau | \mathcal{F}^\tau_t]$$
where $M_{\tau^n} = Y^n_{\tau^n} + \int_0^{\tau^n} f(Y^n_s,Z^n_s)dA^n_s$ and $M_\tau = Y_\tau + \int_0^{\tau} f(Y_s,Z_s)ds $. 
As in Lemma \ref{lemme_cvM}, we prove that ${\displaystyle \forall L, \ \sup_{t \in [0,L]} |M^n_t-M_t| \xrightarrow{\mathbb{P}} 0.}$
At last, arguing like in Lemma \ref{cvY}, $\forall L$, ${\displaystyle \sup_{t \in [0,L]} |Y^n_t-Y_t| \xrightarrow{\mathbb{P}} 0}$. 
In particular, $M^n_0=Y^n_0 \xrightarrow{\mathbb{P}} Y_0=M_0$. Finally, 
$$\forall L, \ \sup_{t \in [0,L]} \left| \int_0^{t \wedge \tau^n} Z^n_s dW^n_s - \int_0^{t \wedge \tau} Z_s dW_s \right| \xrightarrow{\mathbb{P}} 0.$$
Theorem \ref{theo_Donsker} is proved.


\subsection{An example}
\label{Donsker_ex}

Before studying an example with particular stopping times, let us show a technical lemma that will be useful.

\begin{lemme}
\label{hb_tech}
Let $(a^n)_n$ be a sequence of real numbers that converges to $a$. 
We consider the functions $f$ and $f^n$ defined from $\mathbb{D}$ to $\mathbb{R}$ by 
$$f(x) = \inf\{t>0 : x(t) > a\} \text{~~and~~} f^n(x) = \inf\{t \in ]0,n] : x(t) > a^n\} \wedge n$$
with the convention $\inf \emptyset = + \infty$. 
Let $y$ be a continuous process such that 
$$\inf\{t>0 : y(t) > a\}=\inf\{t>0 : y(t) \geqslant a\}.$$
Let $(y^n)_n$ be a sequence of functions of $\mathbb{D}$ that converges to $y$ for the uniform topology on every compact set. 
Then $f^n(y^n) \to f(y)$.
\end{lemme}

\demo
For every $n$, we denote $t_n=f^n(y^n)$. We have $y^n(t_n) \geqslant a^n$ or $t_n=n$. \\
Let $t$ be the limit of a subsequence of $(t_n)$ in $\bar{\mathbb{R}}^+$. Let us show that $t=f(y)$. \\
Without loss of generality, instead of extracting a subsequence, we suppose that $t_n \to t$. As $y$ is continuous and $(y^n)$ converges
uniformly to $y$, $y^n(t_n) \to y(t)$. 
Then, when $n$ tends to $\infty$ in the inequality $y^n(t_n) \geqslant a^n$, we have $y(t) \geqslant a$. So, $t \geqslant f(y)$ because we
have assumed that $\inf\{t>0 : y(t) > a\}=\inf\{t>0 : y(t) \geqslant a\}$. \\
Let us suppose that $t > f(y)$. 
Let us fix $0 < \varepsilon < \frac{t - f(y)}{2}$. 
By definition of $f(y)$, we can find $t_0 \in [f(y),f(y)+\varepsilon[$ such that $y(t_0) >a$. 
We take $\alpha = \frac{y(t_0)-a}{2} \wedge \varepsilon$. Since 
$y^n \to y$, $t_n \to t$ and $a^n \to a$, there exists $n_0$ such that for every $n \geqslant n_0$, $\sup_t|y^n_t - y_t| < \alpha/4$, 
$|t_n - t| < \alpha/4$ and $|a^n-a| < \alpha$. In particular, for every $n \geqslant n_0$, $|y^n(t_0)-y(t_0)| < \alpha/4$. Then,
$$t_n - t_0 > t - \alpha/4 - (f(y) + \varepsilon + \alpha/4) > t -f(y) -\frac{t-f(y)}{2} >0,$$ 
so $t_n > t_0$. \\
On the other hand, $|y^n(t_0)-y(t_0)| < \alpha/4$ and $|a^n-a| < \alpha$, so 
$$y^n(t_0) - a^n > y(t_0) - \alpha/4 - a - \alpha > \frac{3\alpha}{4} >0.$$ 
Hence, $y^n(t_0) > a^n$. 
This is in contradiction with the definition of $t_n$.
So, $t \leqslant f(y)$. \\
Hence, $t=f(y)$ and $f(y)$ is the only possible limit for the convergent subsequences of $(f^n(y^n))_n$. 
So $f^n(y^n) \to f(y)$. 
\findem

\begin{prop}
\label{hb_ex1}
With the notations of Theorem \ref{theo_Donsker}, we assume that $\forall n, \forall k, \varepsilon^n_k=\varepsilon_k$, $\xi = g(W)$ and 
$\xi^n=g(W^n)$ with $g$ bounded and continuous. Let $(a^n)_n$ be a sequence of real numbers that converges to $a$. 
We define the stopping times $(\tau^n)_n$ and $\tau$ as follows : 
$$\tau^n = \inf\{t \in ]0,n] : |W^n_t| > a^n\} \wedge n \text{~~and~~} \tau = \inf\{t>0 : |W_t| > a\}.$$
Then $\left(Y^n_{. \wedge \tau^n}, \int_0^{. \wedge \tau^n} Z^n_s dW^n_s \right)$ converges in law to 
$\left(Y_{. \wedge \tau}, \int_0^{. \wedge \tau} Z_s dW_s \right)$ for the Skorokhod topology. 
\end{prop}

\demo
According to Corollary \ref{hb_cvl}, it is sufficient to prove that $(W^n,\tau^n) \xrightarrow{\mathcal{L}} (W,\tau)$ and that 
$\sup_n \left[ (\tau^n)^2 \right] ^{1/2} < + \infty$. \\

According to Donsker's Theorem, we have the convergence in law of $W^n$ to $W$. Using the Skorokhod representation Theorem, we can 
find a probabilistic space $(\tilde{\Omega}, \tilde{\mathcal{A}}, \tilde{\mathbb{P}})$, processes $\tilde{W}^n$ and $\tilde{W}$ such that 
$\tilde{W}^n \sim W^n$, $\tilde{W} \sim W$ and $\tilde{W}^n \xrightarrow{a.s.} \tilde{W}$. \\
Denoting by $E$ the set $\{\omega : \inf\{t>0:|\tilde{W}_t(\omega)|>a\}\not=\inf\{t>0:|\tilde{W}_t(\omega)|\geqslant a\}\}$, it is quickly 
proved that $\tilde{\mathbb{P}}[E]=0.$
Then, for every $\omega \notin E$, $t \mapsto \tilde{W}_t(\omega)$ is continuous, using Lemma \ref{hb_tech}, for every 
$\omega \notin E$, $f^n(\tilde{W}^n(\omega)) \to f(\tilde{W}(\omega))$. Then, $f^n(\tilde{W}^n) \xrightarrow{a.s.} f(\tilde{W})$. We also have 
$$(\tilde{W}^n,f^n(\tilde{W}^n)) \xrightarrow{a.s.} (\tilde{W},f(\tilde{W})).$$
By construction, $(\tilde{W}^n,f^n(\tilde{W}^n)) \sim (W^n,f^n(W^n))$. Then, we can find a process $Y$ such that 
$(\tilde{W},f(\tilde{W})) \sim (W,Y)$ and $(W^n,f^n(W^n)) \xrightarrow{\mathcal{L}} (W,Y)$. But, by construction, we also have 
$(\tilde{W},f(\tilde{W})) \sim (W,f(W))$. Then $Y = f(W)$ $a.s$. So we have
$$(W^n,f^n(W^n)) \xrightarrow{\mathcal{L}} (W,f(W)), \ i.e. \ (W^n,\tau^n) \xrightarrow{\mathcal{L}} (W,\tau).$$ 

On the other hand, following the lines of the proofs of Proposition 1.16 and Theorem 1.17 in Chung and Zhao \cite{CZ}, we have 
$\sup_n \left[ (\tau^n)^2 \right] ^{1/2} < + \infty.$

The result follows using Corollary \ref{hb_cvl}.
\findem


\section{Stability of BSDEs when the Brownian motion is approximated by a sequence of martingales}
\label{stab_mart}


\subsection{Statement of the problem}
\label{problem}

Let $W$ be a Brownian motion and $\mathcal{F}$ its natural filtration. Let $\tau$ be a $\mathcal{F}$-stopping time almost surely finite. \\
We consider the following BSDE :
\begin{equation}
\label{Z}
Y_{t \wedge \tau} = \xi + \int_{t \wedge \tau}^{\tau} f(r,Y_r,Z_r)dr - \int_{t \wedge \tau}^{\tau} Z_rdW_r,\ t \geqslant 0,
\end{equation}
where $\xi$ is a bounded random variable $\mathcal{F}_\tau$-measurable and for every $(y,z)$, $\{f(t,y,z)\}_t$ is progressively measurable. \\

We approximate this equation on the following way. 
Let $(W^n)_n$ be a sequence of c\`adl\`ag processes and $(\mathcal{F}^n)_n$ the natural filtrations for these processes. We suppose that $(W^n)$ is
a sequence of square integrable $(\mathcal{F}^n)$-martingales which converges in probability to $W$. We don't suppose that $W^n$ has the predictable
representation property. Let $(\tau^n)_n$ be a sequence of $(\mathcal{F}^n)$-stopping times that converges almost surely to $\tau$. \\
Then, we consider the following BSDE :
\begin{equation}
\label{Zn}
Y^n_{t} = \xi^n + \int_{t \wedge \tau^n}^{\tau^n} f^n(r,Y^n_{r-},Z^n_r)d<W^n>_r 
	- \int_{t \wedge \tau^n}^{\tau^n} Z^n_rdW^n_r - (N^n_{\tau^n}-N^n_{t \wedge \tau^n}), t \geqslant 0
\end{equation}
where $(\xi^n)_n$ is a sequence of random variables $(\mathcal{F}^n_{\tau^n})$-measurable, $(N^n)$ is a sequence of $(\mathcal{F}^n)$ martingales orthogonal 
to $(W^{n,\tau^n})$ and for every $(y,z)$, $\{f^n(t,y,z)\}_t$ is progressively measurable with respect to $(\mathcal{F}^n)$. \\

We denote by $\mathcal{S}^p_L$ the set of c\`adl\`ag processes $X$ indexed by $\mathbb{R}^+$ and taking values in $\mathbb{R}$ such that 
$$\|X\|_{\mathcal{S}^p_L} = \mathbb{E} \left[ \sup_{t \in [0,L]} |X_t|^p \right] < + \infty.$$

We put the following assumptions on the martingales and on the terminal conditions :
\begin{center}
\begin{tabular}{ll}
\bf{(H1)} & $(i)$ $\forall L$, $W^n \xrightarrow{\mathcal{S}^2_L} W$, \\
& $(ii)$ $<W^n>_t - <W^n>_s \leqslant \rho(t-s) + a_n$ \\
& ~~~~~~~~~~where $\rho : \mathbb{R}^+ \to \mathbb{R}^+$ with $\rho(0^+)=0$ and $(a_n) \downarrow 0$,\\
& $(iii)$ $\sup_n <W^n>_{\tau^n} < + \infty$. \\
\bf{(H2)} & $(i)$ $\xi^n \xrightarrow{L^2} \xi$, \\
& $(ii)$ $\|\xi\|_\infty + \sup_n \mathbb{E}[|\xi^n|] < \infty$.
\end{tabular}
\end{center}


\subsection{Existence and uniqueness of the solutions for the studied BSDEs}
\label{exist+unic}

Let us begin with the case of the equation (\ref{Z}). \\
Let us put some assumptions on the generator $f$ :
\begin{center}
\begin{tabular}{ll}
\bf{(Hf)} & $(i)$ $f$ is $K$-Lipschitz in $y$ and $z$, \\
& $(ii)$ $f$ is monotone in $y$ in the following way : there exists $\mu > 0$ such that \\
& $\quad \forall (t,y,z),(t,y',z) \in \mathbb{R}^+ \times \mathbb{R}^2, (y-y')(f(t,y,z)-f(t,y',z)) \leqslant - \mu (y-y')^2.$ \\
& $(iii)$ $f$ is bounded.
\end{tabular}
\end{center}

Under these assumptions, according to Theorem 2.1 of Royer in \cite{article_Manuela}, the BSDE (\ref{Z}) has a unique solution $(Y,Z)$ 
(in the sense of Definition \ref{hb_def_sol}) in the set of processes such that $Y$ is continuous and uniformly bounded. \\

Now, let us deal with the equation (\ref{Zn}). \\
Let us put some assumptions on the generators $(f^n)$ :
\begin{center}
\begin{tabular}{ll}
\bf{(Hfn)} & $(i)$ for every $n$, $f^n$ is $K$-Lipschitz in $y$ and $z$, \\
& $(ii)$ $\sup_n \|f^n\| < \infty$.
\end{tabular}
\end{center}

Let us first introduce some notations :\\
$\mathcal{S}^{2,n}$ is the set of processes $Y$ progressively measurable with respect to $\mathcal{F}^{n,\tau^n}$ such that
	$\mathbb{E} \left[ \underset{t \geqslant 0}{\sup} |Y_{t \wedge \tau^n}|^2 \right] < \infty$,\\
$\mathcal{M}^{2,n}$ is the set of predictable processes $Z$ measurable with respect to $\mathcal{F}^{n,\tau^n}$ such that
	$\mathbb{E} \left[ \int_0^{\tau^n} |Z_r|^2 d<W^n>_r \right] < +\infty$, \\
$\mathcal{H}_0^{2,n}$ is the set of squared integrable $\mathcal{F}^{n,\tau^n}$-martingales $M$ such that $M_0=0$. \\

Now, a fixed point argument and an estimation a priori ~like Briand, Delyon and M\'emin in the proof of Theorem 9 in \cite{BDM_mart} give 
the result of existence and uniqueness of the solution of the BSDE (\ref{Zn}) : 

\begin{theo}
Under the assumptions (H1), (H2), (Hf) and (Hfn), the BSDE (\ref{Zn}) has, for n large enough, a unique solution  
$(Y^n_{. \wedge \tau^n},Z^n_{. \wedge \tau^n},N^n_{. \wedge \tau^n})$ in
$\mathcal{S}^{2,n} \times \mathcal{M}^{2,n} \times \mathcal{H}_0^{2,n}$.
\end{theo}


\subsection{Convergence of the solutions}
\label{hb_section_cvsol2}

First, we show a result of stability for the decompositions of the terminal conditions $\xi$ and $\xi^n$. This result will be 
the main argument in the proof of Theorem \ref{theo_mart} about the convergence of the solutions. 

\begin{theo}
\label{hbth5}
We suppose that the conditions (H1) and (H2) are filled. We consider the orthogonal decomposition of  $\xi^n$ with respect to 
$W^n_{. \wedge \tau^n}$, ie $Z^n$ is a predictable $\mathcal{F}^{n,\tau^n}$-measurable process, $N^n$ is a $\mathcal{F}^{n,\tau^n}$-martingale 
orthogonal to $W^n_{. \wedge \tau^n}$ and 
$$M^n_{t \wedge \tau^n} = \mathbb{E}[\xi^n | \mathcal{F}^{n,\tau^n}_t] = M^n_0 + \int_0^{t \wedge \tau^n} Z^n_r dW^n_r + N^n_{t}
= M^n_0 + \int_0^{t \wedge \tau^n} Z^n_r dW^n_r + N^n_{t \wedge \tau^n}.$$
We also consider the representation of $\xi$ as a stochastic integral :
$$M_{t \wedge \tau} = \mathbb{E}[\xi | \mathcal{F}^\tau_t] = M_0 + \int_0^{t \wedge \tau} Z_r dW_r.$$
Then we have the following convergences : for every $L$, 
\begin{displaymath}
\begin{array}{c}
(M^n_{. \wedge \tau^n}, \int_0^{. \wedge \tau^n} Z^n_r dW^n_r, N^n_{. \wedge \tau^n}) 
	\xrightarrow{\mathcal{S}^2_L} (M_{. \wedge \tau}, \int_0^{. \wedge \tau} Z_r dW_r, 0),\\
(\int_0^{. \wedge \tau^n} Z^n_r d<W^n>_r, \int_0^{. \wedge \tau^n} |Z^n_r|^2 d<W^n>_r) 
	\xrightarrow{\mathcal{S}^2_L \times \mathcal{S}^1_L} (\int_0^{. \wedge \tau} Z_r dr, \int_0^{. \wedge \tau} |Z_r|^2 dr).
\end{array}
\end{displaymath}
\end{theo}

\demo
The only noticeable difference with Briand, Delyon and M\'emin's proof of Theorem 5 in \cite{BDM_mart} is that we have to prove the 
convergence of $M^n_{. \wedge \tau^n}$ to $M_{. \wedge \tau}$. \\
According to Proposition 3 in Briand, Delyon and M\'emin in \cite{BDM_mart}, we have the convergence of filtrations 
$\mathcal{F}^n \xrightarrow{w} \mathcal{F}$. As $W$ is continuous and $\tau^n \xrightarrow{\mathbb{P}} \tau$, according to Corollary \ref{cv_proc_arrete_continu}, 
$(\mathcal{F}^{n,\tau^n})$ converges to $\mathcal{F}^\tau$. Moreover, according to (H2), $(\xi^n)$ converges in $L^2$ so in $L^1$ to $\xi$. Then, 
$M^n_{. \wedge \tau^n} \xrightarrow{\mathbb{P}} M_{. \wedge \tau}$, according to Remark 1.2 in Coquet, M\'emin and S\l omi\'nski \cite{cvfiltration}. 
Moreover, the limit is continuous and the sequence is uniformly integrable, so :
$$\forall L, \ M^n_{. \wedge \tau^n} \xrightarrow{\mathcal{S}^2_L} M_{. \wedge \tau}.$$ 
\findem

We are now going to be interested in the convergence of the solutions of the BSDE. \\

Let us do some assumptions of convergence on the generators $f$ and $f^n$.

\begin{center}
\begin{tabular}{ll}
\bf{(H3)} & $\forall (y,z)$, $\{f^n(t,y,z)\}_t$ has c\`adl\`ag trajectories and $f^n(.,y,z) \xrightarrow{\mathcal{S}^2_L} f(.,y,z)$,
$\forall L$.\\
\end{tabular}
\end{center}

\medskip

Now, the following result of convergence of the solutions can be proven in a similar way as in Briand, Delyon and M\'emin in the proof 
of Theorem 12 in \cite{BDM_mart}.
\begin{theo}
\label{theo_mart}
We suppose that (H1), (H2), (H3), (Hf) and (Hfn) are filled. We denote by $(Y^n, Z^n, N^n)$ the solution of the equation (\ref{Zn}) and by 
$(Y, Z)$ these of the equation (\ref{Z}). Then, $(Y^n, Z^n, N^n) \to (Y, Z, 0)$, $i.e.$ $\forall L$, 
\begin{displaymath}
\begin{array}{c}
(Y^n_{. \wedge \tau^n}, \int_0^{. \wedge \tau^n} Z^n_r dW^n_r, N^n_{. \wedge \tau^n}) 
	\xrightarrow{\mathcal{S}^2_L} (Y_{. \wedge \tau}, \int_0^{. \wedge \tau} Z_r dW_r, 0),\\
(\int_0^{. \wedge \tau^n} Z^n_r d<W^n>_r, \int_0^{. \wedge \tau^n} |Z^n_r|^2 d<W^n>_r) 
	\xrightarrow{\mathcal{S}^2_L \times \mathcal{S}^1_L} (\int_0^{. \wedge \tau} Z_r dr, \int_0^{. \wedge \tau} |Z_r|^2 dr).
\end{array}
\end{displaymath}
\end{theo} 

\medskip


\subsection{Application to discretizations}

In this section, we are interested in the case of the approximation of a Brownian motion $W$ by its discretizations $W^n$. More precisely,
we consider an increasing sequence $(\pi^n=\{t_k^n\})_n$ of subdivisions of $\mathbb{R}^+$ with mesh going to 0 and the discretized processes 
$W^n$ are defined by $W^n_{t \wedge \tau^n}=W_{t_k^n}$ if $t_k^n \leqslant t \wedge \tau^n < t_{k+1}^n$. \\

Let $(a^n)_n$ be a sequence of real numbers which decreases to a real number $a$. 
Then, we consider the following random variables : \\
$$\tau = \inf\{t>0 : |W_t| >a\} \text{~~and~~} \tau^n = \inf\{t \in ]0,n] : |W^n_t| >a^n\} \wedge n.$$

As previously, we denote by $\mathcal{F}$ the natural filtration for $W$ and $\mathcal{F}^n$ the natural filtrations of the processes $W^n$. It is clear that
$(W^n)$ is a sequence of $(\mathcal{F}^n)$-martingales, $\tau$ is a $\mathcal{F}$-stopping time and $(\tau^n)_n$ is a sequence of $(\mathcal{F}^n)$-stopping times. \\

Let us show now a result of convergence on these stopping times :
\begin{lemme}
\label{hb_stoptime}
$\tau^n \xrightarrow{a.s.} \tau$.
\end{lemme}

\demo
The proof will be done in two steps. In a first step, we prove that the limit of $(\tau^n)$ is a stopping time. In a second step, 
we shall identify the limit using properties of discretizations. \\

Note that $(\tau^n)_n$ is an almost surely nonincreasing sequence of $(\mathcal{F}^n)$-stopping times lower-bounded by 0. So, $(\tau^n)$ converges 
almost surely to a random variable $\tilde{\tau}$. $\tilde{\tau}$ is a $\mathcal{F}$-stopping time according to the following proposition~:
\begin{prop}
We suppose that, for every $n$, $\mathcal{F}^n \subset \mathcal{F}$ and that the filtration $\mathcal{F}$ is right continuous and complete. 
Let $(\tau^n)_n$ be a sequence of $(\mathcal{F}^n)$-stopping times that converges almost surely to a random variable $\tau$. 
Then $\tau$ is a $\mathcal{F}$-stopping time.
\end{prop}
\demo
Let us fix $t \in \{s : \mathbb{P}[\tau=s]=0\}$. \\
Since $\tau^n \xrightarrow{a.s.} \tau$ and $\mathbb{P}[\tau=t]=0$, we know that $\mathbf{1}_{\{\tau^n \leqslant t\}} \xrightarrow{L^1} \mathbf{1}_{\{\tau \leqslant t\}}$. Then, 
$\mathbb{E}[\mathbf{1}_{\{\tau^n \leqslant t\}}|\mathcal{F}_t] \xrightarrow{L^1} \mathbb{E}[\mathbf{1}_{\{\tau \leqslant t\}}|\mathcal{F}_t].$
But, $\mathbf{1}_{\{\tau^n \leqslant t\}}=\mathbb{E}[\mathbf{1}_{\{\tau^n \leqslant t\}}|\mathcal{F}_t]$ because $\mathcal{F}^n \subset \mathcal{F}$ and $\tau^n$ is a $\mathcal{F}^n$-stopping time. \\
So, by uniqueness of the limit, $\mathbb{E}[\mathbf{1}_{\{\tau \leqslant t\}}|\mathcal{F}_t]=\mathbf{1}_{\{\tau \leqslant t\}}$ $a.s$. Then, $\{\tau \leqslant t\} \in \mathcal{F}_t$. \\
\indent As $\tau$ is a random variable, $\{t : \mathbb{P}[\tau=t] \not= 0 \}$ is countable. 
Let us fix $t$ such that $\mathbb{P}[\tau=t] \not= 0$. We can find a sequence $(t^n)_n$ that decreases to $t$ such that for every $n$, 
$\mathbb{P}[\tau=t^n]=0$. Then $\{\tau \leqslant t\} = \bigcap_n \{\tau \leqslant t^n\}.$
But, for every $n$, $\{\tau \leqslant t^n\} \in \mathcal{F}_{t^n}$.
So $\{\tau \leqslant t\} \in \bigcap_n \mathcal{F}_{t^n}$.
But $\bigcap_n \mathcal{F}_{t^n} = \mathcal{F}_{t^+} = \mathcal{F}_t$ because $\mathcal{F}$ is right continuous. \\
So, for every $t$, $\{\tau \leqslant t\} \in \mathcal{F}_t$, $ie$ $\tau$ is a $\mathcal{F}$-stopping time.
\findem

Using $W^n \xrightarrow{a.s.} W$, $\tau^n \xrightarrow{a.s.} \tilde{\tau}$ and $W$ is continuous, we get $W^n_{\tau^n} \xrightarrow{a.s.} W_{\tilde{\tau}}$. Moreover,
by construction of $\tau^n$, either $\tau^n = n$, or $W^n_{\tau^n} \geqslant a^n$. 
When $n$ tends to $\infty$, either $\tilde{\tau}=+\infty$, or $W_{\tilde{\tau}} \geqslant a$ $a.s$. 
Then, $\tilde{\tau} \geqslant \tau$ $a.s$.\\

Let us fix $\omega$ such that $\tau^n(\omega) \to \tilde{\tau}(\omega)$. 
We suppose that $\tilde{\tau}(\omega) \not= \tau(\omega)$, ie $\tilde{\tau}(\omega) > \tau(\omega)$. 
Then, we can find $t_0 < \tilde{\tau}(\omega)$ such that $W_{t_0} (\omega) \geqslant a$ and $W_{t_0+} (\omega) > a$. As $W$ is right
continuous, we can find $0< \eta <\tilde{\tau}(\omega) - t_0$ such that for every $t \in ]t_0,t_0+\eta[$, $W_t(\omega) > a$.
Let us fix $t_1 \in ]t_0, t_0 +\eta[$. $(a^n)_n$ decreases to $a$, so we can find $n_0$ such that for every $n \geqslant n_0$, $a^n
\leqslant \frac{W_{t_1}(\omega) - a}{2}$. As $W$ is right continuous at time $t_1$, there exists $0< \eta_1 <\tilde{\tau}(\omega) - t_0$ 
such that for every $t \in ]t_1,t_1+\eta_1[$, $W_t(\omega) > \frac{W_{t_1}(\omega) - a}{2}$.
$|\pi^n| \to 0$ so there exists $n_1 \geqslant n_0$ such that for every $n \geqslant n_1$, we can find $t_n \in \pi^n$, 
$t_n \in ]t_1,t_1+\eta_1[$. Then, for every $n \geqslant n_1$, $\tau^n(\omega) \leqslant t_n < \tilde{\tau}(\omega)$. This is in 
contradiction with the fact that $(\tau^n(\omega))_n$ decreases to $\tilde{\tau}(\omega)$. 
So $\tilde{\tau}(\omega) = \tau(\omega)$. Finally, $\tau^n \xrightarrow{a.s.} \tau$. \\
Lemma \ref{hb_stoptime} is proved.
\findem

On the other hand, for every $L$, $\sup_{t \in [0,L]} |W^n_t-W_t| \xrightarrow{a.s.} 0$. Moreover, all the
processes are bounded on $[0,L]$. So the convergence is in $\mathcal{S}^2_L$ (cf section \ref{problem}), 
ie $W^n \xrightarrow{\mathcal{S}^2_L} W$. 
Then, we remark that $<W^n>$ is the discretized process of $<W>$. Let us fix $s \leqslant t$. We can find $i,j \in \mathbb{N}$ such that 
$s \in [t_i^n,t_{i+1}^n[$ and $t \in [t_j^n,t_{j+1}^n[$. Then, $<W^n>_t-<W^n>_s = t_j^n - t_i^n \leqslant t-s + |\pi_n|$ where $|\pi_n|$,
the mesh of the subdivision, goes to 0. 
At last, $<W^n>_{\tau^n} \leqslant \tau^n + |\pi^n|$. Then, 
$$\sup_n <W^n>_{\tau^n} \leqslant \sup_n \tau^n + \sup_n |\pi^n|.$$
It is well known (see e.g. ~Theorem 1.17 in Chung and Zhao \cite{CZ}) that $\mathbb{E}[\tau] < + \infty$. So $\tau$ is almost surely finite. 
As $\tau^n \xrightarrow{a.s.} \tau$, $\sup_n \tau^n < + \infty$ $a.s$. As last, as $(|\pi^n|)_n$ decreases to 0, $\sup_n |\pi^n| < + \infty$ $a.s$. 
Finally, $\sup_n <W^n>_{\tau^n} < + \infty$ $a.s$. \\
So, the assumption (H1) is satisfied. \\

Then, we consider terminal conditions $\xi$ and $(\xi^n)$ and some generators $f$ and $(f^n)$ such that the conditions
(Hf), (Hfn), (H2) and (H3) are filled. \\

Let $(Y,Z)$ be the solution of the BSDE 
$$Y_{t \wedge \tau} = \xi + \int_{t \wedge \tau}^{\tau} f(r,Y_r,Z_r)dr - \int_{t \wedge \tau}^{\tau} Z_rdW_r,\ t \geqslant 0, $$
and $(Y^n,Z^n,N^n)$ the solution of the BSDE
$$Y^n_{t} = \xi^n + \int_{t \wedge \tau^n}^{\tau^n} f^n(r,Y^n_{r-},Z^n_r)d<W^n>_r 
	- \int_{t \wedge \tau^n}^{\tau^n} Z^n_rdW^n_r - (N^n_{\tau^n}-N^n_{t \wedge \tau^n}), t \geqslant 0.$$
These solutions exist and are unique in specified spaces as it was proved in section \ref{exist+unic}. \\

The assumptions of Theorem \ref{theo_mart} are satisfied. So, we have the following convergences : for every $L$, 
\begin{eqnarray*}
\left(Y^n_{. \wedge \tau^n}, \int_0^{. \wedge \tau^n} Z^n_r dW^n_r, N^n_{. \wedge \tau^n} \right) 
	& \xrightarrow{\mathcal{S}^2_L} & \left(Y_{. \wedge \tau}, \int_0^{. \wedge \tau} Z_r dW_r, 0 \right),\\
\left(\int_0^{. \wedge \tau^n} Z^n_r d<W^n>_r, \int_0^{. \wedge \tau^n} |Z^n_r|^2 d<W^n>_r\right) 
	& \xrightarrow{\mathcal{S}^2_L \times \mathcal{S}^1_L} & 
	\left(\int_0^{. \wedge \tau} Z_r dr, \int_0^{. \wedge \tau} |Z_r|^2 dr\right).
\end{eqnarray*}

\noindent $<W^{n,\tau^n}>=<W^n>^{\tau^n}$ according to the following lemma :
\begin{lemme}
Let $M$ be a $\mathcal{F}$-martingale and $\tau$ a $\mathcal{F}$-stopping time. Then we have the following equality : $<M>^\tau = <M^\tau>$.
\end{lemme}
Note that $<W^n>$ is the discretization of $<W>$ : 
$<W^n>_t=t_i^n$ if $t \in [t^n_i,t^n_{i+1}[$. So $<W^{n,\tau^n}>$ is the discretization of $<W^\tau>$. 
Then, using the same arguments as in Lemma \ref{cvZ}, we have :
$$\mathbb{E}\left[\int_0^{\tau \wedge \tau^n} \left|Z^n_{t \wedge \tau^n} - Z_{t \wedge \tau}\right|^2 dr\right] \xrightarrow[n \to + \infty]{} 0.$$
Finally, we have the following convergence for the solutions :
$$\forall L, \ \mathbb{E}\left[ \sup_{t \in [0,L]} |Y^n_{t \wedge \tau^n} - Y_{t \wedge \tau}|^2 
	+ \int_0^{\tau \wedge \tau^n} |Z^n_{t \wedge \tau^n} - Z_{t \wedge \tau}|^2 dr 
	+ \sup_{t \in [0,L]} |N^n_{t \wedge \tau^n}|^2 \right] \xrightarrow[n \to + \infty]{} 0.$$

We have just proved the following theorem :
\begin{theo}
Let $(\pi^n=\{t_k^n\})_n$ be an increasing sequence of subdivisions of $\mathbb{R}^+$ with mesh going to 0 and $W^n$ the discretized associated 
processes of $W$. Let $(a^n)_n$ be a sequence of real numbers which decreases to a real number $a$. 
We consider the stopping times 
$$\tau = \inf\{t>0 : |W_t| >a\}  \text{~~and~~} \tau^n = \inf\{t \in ]0,n] : |W^n_t| >a^n\} \wedge n.$$
Let $(Y,Z)$ be the solution of the BSDE 
$$Y_{t} = \xi + \int_{t \wedge \tau}^{\tau} f(r,Y_r,Z_r)dr - \int_{t \wedge \tau}^{s \wedge \tau} Z_rdW_r,\ t \geqslant 0, $$
and $(Y^n,Z^n,N^n)$ the solution of the BSDE
$$Y^n_{t} = \xi^n + \int_{t \wedge \tau^n}^{\tau^n} f^n(r,Y^n_{r-},Z^n_r)d<W^n>_r 
	- \int_{t \wedge \tau^n}^{\tau^n} Z^n_rdW^n_r - (N^n_{\tau^n}-N^n_{t \wedge \tau^n}), t \geqslant 0.$$
We assume that the conditions (Hf), (Hfn), (H2) and (H3) are satisfied. Then, we have the following convergence for the solutions :
$$\forall L, \ \mathbb{E}\left[ \sup_{t \in [0,L]} |Y^n_{t \wedge \tau^n} - Y_{t \wedge \tau}|^2 
	+ \int_0^{\tau \wedge \tau^n} |Z^n_{t \wedge \tau^n} - Z_{t \wedge \tau}|^2 dr 
	+ \sup_{t \in [0,L]} |N^n_{t \wedge \tau^n}|^2 \right] \xrightarrow[n \to + \infty]{} 0.$$
\end{theo}  


\appendix
\section{About stopped filtrations and stopped processes}
\label{cv_stop_filtr}

In their paper \cite{caractF_tau}, Haezendonck and Delbaen give the following characterization of the $\sigma$-field $\mathcal{F}_\tau$ when $\mathcal{F}$ is
the natural filtration of a process $X$ : 

\begin{prop}
\label{caractF_tau}
Let $X$ be a c\`adl\`ag process, $\mathcal{F}$ the natural filtration of $X$ and $\tau$ a $\mathcal{F}$-stopping time. 
Then $\mathcal{F}_\tau =\sigma(\{X_{\tau \wedge s}, s \geqslant 0 \})$. 
\end{prop}

This characterization shows that, if $\mathcal{F}$ is the natural filtration of $X$ and $\tau$ a $\mathcal{F}$-stopping time, the stopped filtration 
$\mathcal{F}^\tau$ is the natural filtration of the stopped process $X^\tau$. \\

The notions of convergence of filtrations and of $\sigma$-fields have been firstly defined in Hoover \cite{Hoover} and then in a slightly
different way in Coquet, M\'emin and S\l ominski \cite{cvfiltration}. In \cite{cvfiltration}, the filtrations are indexed by  a finite 
interval time $[0,T]$. We generalise it to the case of filtrations indexed by $\mathbb{R}^+$.  

\begin{definition}
We say that $(\mathcal{F}^n)$ converges to $\mathcal{F}$ if for every $A \in \mathcal{F}_\infty$, the sequence of processes $(\mathbb{E}[\mathbf{1}_A | \mathcal{F}^n_.])_n$ 
converges in probability to $\mathbb{E}[\mathbf{1}_A | \mathcal{F}_.]$ for the Skorokhod topology. We denote $\mathcal{F}^n \xrightarrow{w} \mathcal{F}$.
\end{definition}

\begin{definition}
We say that the sequence of $\sigma$-fields $(\mathcal{B}_n)$ converges to the $\sigma$-field $\mathcal{B}$ if for every $A \in \mathcal{B}$, the sequence of random
variables $(\mathbb{E}[\mathbf{1}_A | \mathcal{B}^n])_n$ converges in probability to $\mathbf{1}_A$. We denote $\mathcal{B}_n \to \mathcal{B}$.
\end{definition}

The following lemma shows that, when holds convergence of filtrations, to get the convergence of associated stopped filtrations
we just have to check the convergence of the terminal $\sigma$-fields. More precisely,

\begin{lemme}
\label{cvtribufiltration}
Let $(\mathcal{F} ^n)$ be a sequence of filtrations that converges to the filtration $\mathcal{F}$. Let $(\tau ^n)$ be a sequence of $\mathcal{F} ^n$-stopping times
that converges in probability to a $\mathcal{F}$-stopping time $\tau$. If the convergence of the $\sigma$-fields $(\mathcal{F}^n_{\tau^n})_n$ to 
$\mathcal{F}_\tau$ holds, then also holds the convergence of the filtrations $(\mathcal{F}^{n,\tau^n})_n$ to $\mathcal{F}^{\tau}$.
\end{lemme}

\demo
Let us fix $B \in \mathcal{F}_\tau$. \\
$\mathcal{F}^n_{\tau^n} \to \mathcal{F}_\tau$, so by definition $\mathbb{E}[\mathbf{1}_B | \mathcal{F}^n_{\tau^n}] \xrightarrow{\mathbb{P}} \mathbf{1}_B$. Moreover, this convergence holds in $L^1$ since the 
sequence is uniformly integrable. As $\mathcal{F}^n \to \mathcal{F}$, according to Remark 1.2 in Coquet, M\'emin and S\l ominski \cite{cvfiltration}, we have
\begin{equation}
\label{eq1}
\mathbb{E}[\mathbb{E}[\mathbf{1}_B | \mathcal{F}^n_{\tau^n}]|\mathcal{F}^n_.] \xrightarrow{\mathbb{P}} \mathbb{E}[\mathbf{1}_B | \mathcal{F}_.] \text{~~for the Skorokhod topology.}
\end{equation}
For every $t$, we have the relations 
$\mathbb{E}[\mathbb{E}[\mathbf{1}_B | \mathcal{F}^n_{\tau^n}]|\mathcal{F}^n_t] = \mathbb{E}[ \mathbf{1}_B | \mathcal{F}^n_{\tau^n \wedge t}] = \mathbb{E}[ \mathbf{1}_B | \mathcal{F}^{n, \tau^n}_t]$ 
and $\mathbb{E}[\mathbf{1}_B | \mathcal{F}_t] = \mathbb{E}[\mathbf{1}_B | \mathcal{F}^{\tau}_t]$ using Proposition 1.2.17 in Karatzas and Shreve \cite{KS}.
So the convergence (\ref{eq1}) can be written on the following way : 
$\mathbb{E}[\mathbf{1}_B | \mathcal{F}^{n,\tau^n}_.] \xrightarrow{\mathbb{P}} \mathbb{E}[\mathbf{1}_B | \mathcal{F}^{\tau}_.]$ for the Skorokhod topology.
Lemma \ref{cvtribufiltration} is proved.
\findem

Using the same arguments as in the proof of Lemma 3 in \cite{cvfiltration}, we have the following characterization of convergence of 
$\sigma$-fields : 

\begin{lemme}
\label{caractcvtribu}
Let $Y$ be a c\`adl\`ag process, $\mathcal{A} = \sigma( \{ Y_t, t \geqslant 0\})$ and $(\mathcal{A}^n)$ be a sequence of $\sigma$-fields.
The following assumptions are equivalent : \\
$i)$ $\mathcal{A}^n \to \mathcal{A}$, \\
$ii)$ $\mathbb{E}[f(Y_{t_1}, \ldots, Y_{t_k}) | \mathcal{A}^n] \xrightarrow{\mathbb{P}} f(Y_{t_1}, \ldots, Y_{t_k})$ for every bounded continuous function 
$f : \mathbb{R}^k \to \mathbb{R}$ and $t_1, \ldots, t_k$ continuity points of $Y$.
\end{lemme}

Then, with the characterization of Proposition \ref{caractF_tau}, we can show a link between convergence of stopped processes and 
convergence of stopped filtrations :

\begin{theo}
\label{cv_proc_arrete}
Let $(X^n)$ and $X$ be c\`adl\`ag processes, $(\mathcal{F}^n)$ and $\mathcal{F}$ their natural filtrations. Let $(\tau^n)$ be a sequence of $(\mathcal{F}^n)$-stopping
times that converges in probability to a $\mathcal{F}$-stopping time $\tau$. We suppose that $X^{n,\tau^n} \xrightarrow{\mathbb{P}} X^\tau$ for the Skorokhod topology
and that $\mathcal{F}^n \xrightarrow{w} \mathcal{F}$. Then $\mathcal{F}^{n,\tau^n} \xrightarrow{w} \mathcal{F}^\tau$.
\end{theo}

\demo
As $\tau$ and $\tau^n$ are respectively $\mathcal{F}$ and $\mathcal{F}^n$-stopping times, according to Proposition \ref{caractF_tau}, we have the equalities 
$\mathcal{F}_\tau =\sigma(\{X_{\tau \wedge s}, s \geqslant 0 \})$ and $\mathcal{F}^n_{\tau^n} =\sigma(\{X^n_{\tau^n \wedge s}, s \geqslant 0 \})$. \\
Let $t_1, \ldots t_k$ be points of continuity of $X^\tau$ and $f:\mathbb{R}^k \to \mathbb{R}$ be a bounded continuous function. \\
As $X^{n,\tau^n} \xrightarrow{\mathbb{P}} X^\tau$, we have :
$(X^n_{t_1 \wedge \tau^n}, \ldots , X^n_{t_k \wedge \tau^n}) \xrightarrow{\mathbb{P}} (X_{t_1 \wedge \tau}, \ldots , X_{t_k \wedge \tau}).$
As $f$ is bounded and continuous, we have :
\begin{equation}
\label{cvf2}
f( X^n_{t_1 \wedge \tau^n}, \ldots , X^n_{t_k \wedge \tau^n}) \xrightarrow{L^1} f( X_{t_1 \wedge \tau}, \ldots , X_{t_k \wedge \tau}).
\end{equation}
Finally,
\begin{eqnarray*}
\lefteqn{ \mathbb{P}[|\mathbb{E}[f( X_{t_1 \wedge \tau}, \ldots , X_{t_k \wedge \tau})|\mathcal{F}^n_{\tau^n}]
		- f( X_{t_1 \wedge \tau}, \ldots , X_{t_k \wedge \tau})| \geqslant \eta] } \\
& \leqslant & 
\mathbb{P}[|\mathbb{E}[f( X_{t_1 \wedge \tau}, \ldots , X_{t_k \wedge \tau})|\mathcal{F}^n_{\tau^n}] 
	- \mathbb{E}[f( X^n_{t_1 \wedge \tau^n}, \ldots , X^n_{t_k \wedge \tau^n})|\mathcal{F}^n_{\tau^n}]| \geqslant \eta/2] \\
& & + \mathbb{P}[|\mathbb{E}[f( X^n_{t_1 \wedge \tau^n}, \ldots , X^n_{t_k \wedge \tau^n})|\mathcal{F}^n_{\tau^n}]
	- f( X_{t_1 \wedge \tau}, \ldots , X_{t_k \wedge \tau})| \geqslant \eta/2] \\
& \leqslant &
\frac{4}{\eta} \mathbb{E}[|f( X^n_{t_1 \wedge \tau^n}, \ldots , X^n_{t_k \wedge \tau^n})
	- f( X_{t_1 \wedge \tau}, \ldots , X_{t_k \wedge \tau})|] \\
& \to & 0 \text{~~according to (\ref{cvf2}).}
\end{eqnarray*}
Then, $\mathcal{F}^n_{\tau^n} \to \mathcal{F}_{\tau}$, according to Lemma \ref{caractcvtribu} and, using Lemma \ref{cvtribufiltration}, 
$\mathcal{F}^{n,\tau^n} \xrightarrow{w} \mathcal{F}^\tau$. \\
Theorem \ref{cv_proc_arrete} is proved
\findem

Let us show a Corollary when the limit is a continuous process.

\begin{cor}
\label{cv_proc_arrete_continu}
Let $(X^n)$ be a sequence of c\`adl\`ag processes and $X$ a continuous process, $(\mathcal{F}^n)$ and $\mathcal{F}$ the associated filtrations. 
Let $(\tau^n)$ be a sequence of $\mathcal{F}^n$-stopping times that converges in probability to a $\mathcal{F}$-stopping time $\tau$. We suppose that 
$X^n \xrightarrow{\mathbb{P}} X$ for the Skorokhod topology and that $\mathcal{F}^n \xrightarrow{w} \mathcal{F}$. 
Then we have the convergence of the stopped filtrations $\mathcal{F}^{n,\tau^n} \xrightarrow{w} \mathcal{F}^\tau$.
\end{cor}

\demo
According to Theorem \ref{cv_proc_arrete}, we just have to prove that $X^{n,\tau^n} \xrightarrow{\mathbb{P}} X^\tau$.\\
By definition of the Skorokhod topology, we have to prove that $\forall L \in \mathbb{N}$, $\sup_{t \in [0,L]} |X^{n,\tau^n}_t - X^\tau_t| \xrightarrow{\mathbb{P}} 0$.
\\
Let us fix $L \in \mathbb{N}$ and $\eta > 0$. We have : \\
\begin{eqnarray*}
\lefteqn{\mathbb{P}\left[\sup_{t \in [0,L]}|X^{n,\tau^n}_t - X^\tau_t| \geqslant \eta\right]} \\
& \leqslant & \mathbb{P}\left[\sup_{t \in [0,L]}|X^n_{t \wedge \tau^n} - X_{t \wedge \tau^n}| \geqslant \eta/3\right] 
	+ \mathbb{P}\left[\sup_{t \in [0,L]}|X_{t \wedge \tau^n} - X_{t \wedge \tau}|\mathbf{1}_{|\tau^n-\tau| < \alpha} \geqslant \eta/3\right]\\
&& \quad \quad \quad \quad 
	+ \mathbb{P}\left[\sup_{t \in [0,L]}|X_{t \wedge \tau^n} - X_{t \wedge \tau}|\mathbf{1}_{|\tau^n-\tau| \geqslant \alpha} \geqslant \eta/3\right] \\
& \leqslant & \mathbb{P}\left[\sup_{t \in [0,L]}|X^n_{t} - X_{t}| \geqslant \eta/3\right]
	+ \varepsilon + \mathbb{P}[|\tau^n-\tau| \geqslant \alpha] \\
& \xrightarrow[n \to +\infty]{} & 0
\end{eqnarray*}
because $X^n \xrightarrow{\mathbb{P}} X$, $X$ is continuous on the compact $[0,L]$ and $\tau^n \xrightarrow{\mathbb{P}} \tau$. 
So $X^{n,\tau^n} \xrightarrow{\mathbb{P}} X^\tau$.
\findem

\end{document}